\documentclass{amsart}
\usepackage{amssymb,latexsym}
\usepackage{amsmath}
\usepackage{epsfig}
\usepackage{graphics}
 \usepackage{lscape} 
\usepackage{enumerate}
\usepackage{amsthm}
\usepackage{hyperref}
\usepackage{color}

\newcommand{\be}{\begin{enumerate}}  \newcommand{\ee}{\end{enumerate}}
\newcommand{\beqt}{\begin{equation}}  \newcommand{\eeqt}{\end{equation}}
\newcommand{\beq}{\begin{eqnarray}}  \newcommand{\eeq}{\end{eqnarray}}
\newcommand{\beQ}{\begin{eqnarray*}} \newcommand{\eeQ}{\end{eqnarray*}}

\newcommand{\M}{\mathbb{M}}

\newcommand{\CC}{{\mathbb C}}

\newcommand{\RR}{{\mathbb R}}

\newcommand{\SSS}{{\mathbb S}}

\newcommand{\Ekt}{\mathbb{E}(\kappa,\tau)}

\newcommand{\lto}{\ensuremath{\longrightarrow}}

\renewcommand{\div}{\mathrm{div}}

\newcommand{\Chi}{\mathfrak{X}}

\newcommand{\id}{\mathrm{Id}}

\newcommand{\lgra}{\longrightarrow}

\newtheorem{example}{Examples}[section]
\newtheorem{thm}{Theorem}[section]
\newtheorem{lemma}[thm]{Lemma}
\newtheorem{prop}[thm]{Proposition}
\newtheorem{cor}[thm]{Corollary}
\newtheorem{remark}[thm]{Remark}
\newtheorem{remarks}[thm]{Remarques}
\newtheorem{definition}[thm]{Definition}
\newtheorem{notation}[thm]{Notation}
\newtheorem{exabout:ample}[thm]{Exemple}

\newcommand{\function}[5]
{\begin{eqnarray*}\begin{array}{r@{}ccl}
 #1\;\colon\;  & #2 &\lto & #3 \\[.05cm]  
  & #4 &\longmapsto  & #5 
\end{array}\end{eqnarray*}
}

\begin{document}

\title{Isometric immersions into manifolds with metallic structures}
\author{Julien Roth and Abhitosh Upadhyay}


\maketitle

\begin{abstract}
We consider submanifolds into Riemannian manifold with metallic structures. We obtain some new results for hypersurfaces in these spaces and we express the fundamental theorem of submanifolds into products spaces in terms of metallic structures. Moreover, we define new structures called {\it complex metallic structures}. We show that these structures are linked with complex structures. Then, we consider submanifolds into Riemannian manifold with such structures  with a focus on invariant submanifolds and hypersurfaces. We also express in particular the fundamental theorem of submanifolds of complex space form in terms of complex metallic structures.
\end{abstract}


\section{Introduction}
A classical problem in submanifold theory consists in determining when a given Riemannian manifold $(M^ng)$ can be immersed (at least locally) into a fixed Riemannian manifold $(\overline{M},\overline{g})$. The Gauss, Codazzi and Ricci equations give the relation between curvatures of the submanifold and the second fundamental forms of the ambient manifold. Conversely, in a large variety of cases, when these three fundamental equations can be written in a intrinsic way, that is in terms of quantity defined only on $M$ and the vector bundle over, possibly with some additional conditions, then it is possible to construct a local isometric immersion into the desired ambient space. The first result in this direction is the classical fundamental theorem of surfaces proven by Bonnet \cite{Bo} which states that if a Riemannian surface $(\Sigma,g)$ endowed with a symmetric tensor $B$ satisfies the Gauss and Codazzi equations, then $\Sigma$ can be isometrically immersed into $\RR^3$ with $B$ as second fundamental form. This result has been generalized later in many cases like for higher dimension and codimension submanifolds in real space forms \cite{Ten}, 3-homogeneous manifolds \cite{Da,Da2}, product spaces \cite{Ko,LR,LTV,Ro2}, warped products \cite{LO} and other ambient spaces with sufficient homogeneity so that the Gauss, Coadzzi and Ricci equations can be expressed in an appropriate manner \cite{PT}.\\ \\
In the present article, we are interested in submanifolds into Riemannian metallic manifold.\linebreak A {\it metallic structure} on manifold $M$ is a $(1,1)$-tensor $J$ over $TM$ satisfies the so-called metallic equation $J^2-pJ-q\id=0$, where $p,q$ are two positive integers. Metallic structures appear as particular case of polynomial structures on manifold introduced by Goldberg and Yano \cite{GY}. If $p=q=1$, then $J$ is called a {\it golden structure} since it is a solution of the well known golden equation $J^2-J-\id=0$. Moreover, if $M$ is endowed with a Riemiannan metric, we say that $J$ is a Riemannian metallic structure if $J$ is symmetric with respect to $g$. Submanifolds in metallic Riemannian manifold, in particular, golden Riemannian manifold has been considered only very recently (see \cite{HC2,HC,OP}). Metallic structures are highly linked to product structures and in the present paper, we will explain how to obtain a fundamental theorem of submanifolds for products of two real space forms in terms of the canonical metallic structure (Theorem \ref{thm1}). Then, as applications, we recover the spinorial version for the fundamental theorem of surfaces into the $4$-dimensional products as well as the existence of associated families of minimal surfaces in terms of metallic structures (Theorem \ref{thm2}).\\
In the second part of the article, we introduce new structures called {\it complex metallic structures} which are satisfying the second degree equation $J^2+aJ+b\id=0$ with $a,b>0$ so that $4b-a^2>0$. We show that such structures are in relation with complex structures. Then, we study submanifolds in Riemannian manifold carrying such structures. We are able to prove comparable results than for metallic structures.  In particular, we prove a fundamental theorem of submanifold, here for complex space forms (Theorem \ref{thm1comp}). We also deduce a spinorial version for invariant submanifolds in the complex projective space $\CC P^2$. \\
For both metallic structures and complex metallic structures, examples are given.
\section{Preliminaries}

\subsection{Riemannian metallic structures}\label{secmetallic}
Let $p,q$ be two positive integers. The positive solution of the equation $x^2-px-q=0$, denoted by $\sigma_{p,q}$, is called the $(p,q)$-metallic mean or $(p,q)$-number. Clearly we have 
$$\sigma_{p,q}=\dfrac{p+\sqrt{p^2+4q}}{2}.$$
The denomination {\it metallic} comes obviously for $p=q=1$, i.e., $\sigma_{1,1}=\frac{1+\sqrt{5}}{2}\phi$ which is the well known golden number related to fibonacci numbers. Moreover, for $p=2,q=1$, $\sigma_{2,1}=1+\sqrt{2}$ which is the silver number appearing in the study of the so-called Pell numbers, etc.. One can refer to \cite{dS} for more details about metallic numbers and their link to number theory or fractal geometry.\\
It is to note that for what we will do, there is no restriction for $p$ and $q$ to be integers, so we will consider that $p$ and $q$ are two positive numbers.\\
Now, let $(M,g)$ be a Riemannian manfiold and $p,q$ be two positive integers. We say that the $(1,1)$-tensor $J$ over $M$ is an almost $(p,q)$-metallic structure if it satisfies the metallic equation $J^2-pJ-q\id=0$. Moreover, $J$ is said to be a Riemaniann $(p,q)$-metallic structure if $J$ is compatible with the metric $g$, that is $g(JX,Y)=g(X,JY)$ for any $X,Y\in \Gamma(TM)$. Further, $J$ is said to be parallel if $J$ is parallel with respect to the Levi-Civita connection $\nabla$. There is a natural connection between metallic structures and product structures since every product structure $F$ induces two metallic structures given by
$$J_1=\frac{p}{2}\id+\frac{2\sigma_{p,q}-p}{2}F\quad\text{and}\quad J_2=\frac{p}{2}\id-\frac{2\sigma_{p,q}-p}{2}F.$$
Conversely, every metallic structure $J$ induces two product structures given by
$$F_{\pm}=\pm\left( \frac{2}{2\sigma_{p,q}-p}J-\frac{p}{2\sigma_{p,q}-p}\id \right).$$
Moreover, the two projections associated with the decomposition of the tangent space by the product structure are given by
$$\pi_1=\frac{\sigma_{p,q}}{2\sigma_{p,q}-p}\id-\frac{1}{2\sigma_{p,q}-p}J\quad\text{and}\quad \pi_2=\frac{\sigma_{p,q}-p}{2\sigma_{p,q}-p}\id+\frac{1}{2\sigma_{p,q}-p}J.$$
For basic examples about metallic structures, one can refer to \cite{CH}. We will also detail some examples at the end of Section \ref{secpq}.
\subsection{Fundamental equations of submanifolds in metallic manifolds}\label{secsubmdfpq}
Now, let us consider a Riemannian manifold $(M^n,g)$ isometrically immersed into a $(n+m)$-dimensional Riemannian manifold $(\widetilde{M},\widetilde{g})$ endowed with a Riemannin parallel $(p,q)$-metallic structure $J$. We denote by $E$ the normal bundle which is equipped with an induced metric $g^E$ and the induced compatible normal connection $\nabla^E$. Then, the metallic structure $J$ induces the existence of four operators $P: TM \longrightarrow TM$, $Q: TM \longrightarrow E$, $R:E \longrightarrow TM$ and $S:E \longrightarrow E$, so that with respect to the decomposition $T\widetilde{M}=TM\oplus E$, $J$ is given over $M$ by
$J=\left( 
\begin{array}{cc}
P&R\\
Q&S
\end{array}\right).$
Then, the operators $P, Q, R$ and $S$ satisfy the following equations.
\begin{prop}\label{allrelations}
For all $X,Y\in TM$ and all $\xi,\nu\in E$, we have
\begin{align}
&P^2+R\circ Q=pP+q\id_{TM},&\label{relation1.1}\\
&Q\circ P+S\circ Q=pQ,&\label{relation1.2}\\
&P\circ R+R\circ S=pR,&\label{relation1.3}\\
&S^2+Q\circ R=pS+q\id_{E},&\label{relation1.4}\\
&g(PX,Y)=g(X,PY),&\label{relation1.5}\\
&g(QX,\xi)=g(X,R\xi),&\label{relation1.6}\\
&g(S\xi,\nu)=g(\xi,S\nu),&\label{relation1.7}\\
&\nabla_X(PY)-P(\nabla_XY)=A_{QY}X+R(B(X,Y)),& \label{relation2.1}\\
&\nabla^{\perp}_X(QY)-Q(\nabla_XY)=S(B(X,Y))-B(X,PY),& \label{relation2.2}\\
&\nabla_X^{\perp}(S\nu)-S(\nabla^{\perp}_X\nu)=-B(R\nu,X)-Q(A_{\nu}X),& \label{relation2.3}\\
&\nabla_X(R\nu)-R(\nabla^{\perp}_X\nu)=-P(A_{\nu}X)+A_{S\nu}X.&\label{relation2.4}
\end{align}
{\it Proof:} Writing $J$ as a matrix by blocks with respect to the decomposition $T\widetilde{M}=TM\oplus E$, we have
$$J=\left( 
\begin{array}{cc}
P&R\\
Q&S
\end{array}\right)$$
and so 
$$J^2=\left( 
\begin{array}{cc}
P^2+R\circ S&P\circ R+ R\circ S\\ \\
Q\circ P+S\circ Q& Q\circ R+S^2
\end{array}
\right).$$
\end{prop}
The identities \eqref{relation1.1}-\eqref{relation1.4} are immediate from this and the relation $J^2=pJ+q\id$. Moreover, the relations \eqref{relation1.5}-\eqref{relation1.7} come directly from $g(JX,Y)=g(X,JY)$ for any $X,Y\in T\widetilde{M}$. 

  Finally, \eqref{relation2.1}-\eqref{relation2.4} are consequences of the fact that $J$ is parallel. Indeed, we have for any \linebreak $X,Y\in TM$ and $\nu\in E$,

\beqt\label{nablaparallel}
0=\left(\overline{\nabla}J\right)(Y+\nu)=\overline{\nabla}_X(JY)-J(\overline{\nabla}_XY)+\overline{\nabla}_X(J\nu)-J(\overline{\nabla}_X\nu).
\eeqt
Now, we recall that for any $X,Z\in TM$ and $\xi\in E$, we have
$$\overline{\nabla}_XZ=\nabla_XZ+B(X,Z)\quad\text{and}\quad \overline{\nabla}_X\nu=\nabla^{\perp}_X\xi-A_{\xi}X.$$
Hence \eqref{nablaparallel} becomes
\beQ
0&=&\overline{\nabla}_X(PY)+\overline{\nabla}_X(QY)-P((\overline{\nabla}_XY)^{\top})-Q((\overline{\nabla}_XY)^{\top})\\
&&+\overline{\nabla}_X(R\nu)+\overline{\nabla}_X(S\nu)-R((\overline{\nabla}_X\nu)^{\perp})-S((\overline{\nabla}_X\nu)^{\perp})\\
&=&\nabla_X(PY)+B(X,PY)+\nabla^{\perp}_X(QY)-A_{QY}X-P(\nabla_XY)-Q(\nabla_XY)\\
&&\nabla_X(R\nu)+B(X,R\nu)+\nabla^{\perp}(S\nu)-A_{S\nu}X-R(\nabla^{\perp}_X\nu)-S(\nabla^{\perp}_X\nu).
\eeQ
For $\nu=0$, the tangential and normal parts of the last equality are relations \eqref{relation2.1} and \eqref{relation2.2} whereas for $Y=0$, they give  \eqref{relation2.3} and \eqref{relation2.4}. This concludes the proof of the proposition.
\hfill$\square$
\begin{remark}
We want to point out that in \cite{HC}, the relations \eqref{relation1.1}-\eqref{relation1.7} are given with a different expression, but are equivalent.
\end{remark}
\begin{remark}
Using \eqref{relation1.1}-\eqref{relation1.4}, we can show easily that the relations \eqref{relation1.5}-\eqref{relation1.7} imply the following three identities:
\begin{enumerate}[(i)]
\item $g(PX,PX)+g(QX,QX)=pg(X,PX)+qg(X,Y)$,
\item $g(R\xi,R\nu)+g(S\xi,S\nu)=pg(X\xi,S\nu)+qg(\xi,\nu)$,
\item $g(PX,R\xi)+g(QX,S\xi)=pg(X,R\xi)+pg(QX,Y)$.
\end{enumerate}
\end{remark}

We finish this section by considering two particular cases, namely the hypersurfaces and the invariant submanifolds.
\begin{definition}
A submanifold $M$ into the Riemannian metallic manifold $(\widetilde{M},J)$ is called invariant with respect to $J$ if $J(T_xM)\subset T_xM$ for all $x\in M$.
\end{definition}
Then, we have the following:
\begin{prop}\label{allrelationsinv}
If $M$ is an invariant submanifold with respect to $J$, then the operators $Q$ and $R$ vanish and the operators $P$ and $S$ satisfy the following equations for all $X,Y\in TM$ and all $\xi,\nu\in E$.
\begin{align}
&P^2=pP+q\id_{TM},&\label{relation1.1inv}\\
&S^2=pS+q\id_{E},&\label{relation1.4inv}\\
&g(PX,Y)=g(X,PY),,&\label{relation1.5inv}\\
&g(S\xi,\nu)=g(\xi,S\nu),,&\label{relation1.7inv}\\
&\nabla_X(PY)-P(\nabla_XY)=0,& \label{relation2.1inv}\\
&S(B(X,Y))=B(X,PY),& \label{relation2.2inv}\\
&\nabla_X^{\perp}(S\nu)-S(\nabla^{\perp}_X\nu)=0,& \label{relation2.3inv}\\
&P(A_{\nu}X)=A_{S\nu}X.&\label{relation2.4inv}
\end{align}
In particular, $P$ and $S$ are metallic structures (for the same metallic equation as $J$), respectively, on $TM$ and $E$.
\end{prop}
{\it Proof:} The proof is immediate from Proposition \ref{allrelations} with the fact that $Q=0$ and $R=0$.\hfill$\square$\\ \\ 
Now, we consider hypersurfaces. In this case, it is more convenient to consider the real-valued second fundamental form by taking the scalar product with the unit normal $\nu$. Therefore, the metallic structure $J$ on $\widetilde{M}$ implies the existence of a field of symmetric operators $P:TM\longrightarrow TM$, a vector field $V\in\Gamma(TM)$ and a smooth function $f$ on $M$. Note that $V$ and correspond to the tensosr $R$ and $S$ respectively in this case. The tensor $Q$ is just the dual $1$-form associated to $V$. These three objects satisfiy the following relations
\begin{prop}\label{allrelationshyp}
For all $X,Y\in TM$ of a hypersurface $M$, $P$, $V$ and $f$ satisfy the following relations:
\begin{align}
&P^2+\langle V,\cdot\rangle V=pP+q\id_{TM},&\label{relation1.1hyp}\\
&PV+fV=pV,&\label{relation1.2hyp}\\
&f^2+\|V\|^2=pf+q,&\label{relation1.4hyp}\\
&g(PX,Y)=g(X,PY),,&\label{relation1.5hyp}\\
&\nabla_X(PY)-P(\nabla_XY)=\langle V,Y\rangle AX+\langle AX,Y\rangle V,& \label{relation2.1hyp}\\
&\nabla_XV=-P(AX)+fAX,& \label{relation2.3hyp}\\
&df(X)=-2\langle AX,V\rangle.&\label{relation2.4hyp}
\end{align}

\end{prop}
{\it Proof:} Here again, the proof is immediate from Proposition \ref{allrelations}.\hfill$\square$\\ \\ 
In \cite{HC}, the authors give necessary and sufficient conditions for non-invariant hypersurfaces to be totally geodesic. Namely, they prove that if $M$ is non-invariant, that is $V\neq0$, then $M$ is totally geodesic if and only if $P$ is parallel or equivalentely $V$ is parallel. Here, we give two results for non-invariant hypersurfaces. Namely, we show that if either $J(\nu)$ is tangent $(\nu$ is the unit normal to the hypersurface) or $J(V)$ is normal then the hypersurface have vanishing Gauss-Kronecker curvature. It is to note that both assumptions, $J(\nu)$ is tangent or $J(V)$ is normal, imply that the hypersurface is non-invariant. Precisely, we have the following two results.
\begin{prop}\label{prop1}
Let $(M^n,g)$ be a hypersurface of a Riemannian metallic manifold $(\widetilde{M},\widetilde{g},J)$ with unit normal vector field $\nu$. Let $P:TM\longrightarrow TM$, $V\in\Gamma(TM)$ and $f\in\mathcal{C}^{\infty}(M)$ induced by $J$. If $V$ has no zero and $J(V)$ is normal, then $A(V)=0$, $\|V\|^2=q$ and the shape operator of the immersion is given by
$$ A(X)=-\frac{P(\nabla_XV)}{\|V\|^2}.$$
In particular, the Gauss-Kronecker curvature of $M$ vanishes identically.
\end{prop}
{\it Proof:} We assume that $J(V)$ is normal, then $P(V)=0$. Hence, from \eqref{relation1.2hyp}, we get $fV=pV$. Since $V$ has no zeros, then the function $f$ is constant equal to $p$. Thus, \eqref{relation1.4hyp} becomes $\|V\|^2=q$. Moreover, since $f$ is constant, from \eqref{relation2.4hyp}, for any $X\in \Gamma(TM)$, we have $\langle A(X),V\rangle=0$ or equivalentely $\langle A(V),X\rangle=0$ by the symmetry of $A$. Hence, $A(V)=0$. Therefore, $A$ has a non-trivial kernel and so, the Gauss -Kronecker curvature of $M$ vanishes identically. Finally, from \eqref{relation2.1hyp} for $X=V$, using $P(V)=0$ and $\|V\|^2=q$ give $A(X)=-\frac{P(\nabla_XV)}{\|V\|^2}$.
\hfill$\square$\\ \\
In the same spirit, we have this second proposition.
\begin{prop}\label{prop2}
Let $(M^n,g)$ be a hypersurface of a Riemannian metallic manifold $(\widetilde{M},\widetilde{g},J)$ with unit normal vector field $\nu$. Let $P:TM\longrightarrow TM$, $V\in\Gamma(TM)$ and $f\in\mathcal{C}^{\infty}(M)$ induced by $J$. If $J(\nu)$ is tangent, then the shape operator of the immersion is given by
$$ A(X)=-\frac{P(\nabla_XV)-p\nabla_XV}{\|V\|^2}.$$
In particular, the Gauss-Kronecker curvature of $M$ vanishes identically.
\end{prop}
{\it Proof:} The proof is completely analogous to the proof of Proposition \ref{prop1} with the only difference that $f=0$ since $J(\nu)$ is tangent. We also get by \eqref{relation1.4hyp} that $\|V\|^2=q$ and from \eqref{relation2.4hyp}, that $A(V)=0$. Hence, we deduce that $P(V)=pV$ and $ A(X)=-\frac{P(\nabla_XV)-p\nabla_XV}{\|V\|^2}$ from \eqref{relation1.2hyp} and \eqref{relation2.1hyp}, respectively.
\hfill$\square$\\ \\
In \cite{HC}, the authors gave some sufficent conditions for non-invariant hypersurfaces to be minimal. We can prove the following necessary and sufficient condition, for both cases, if $J(V)$ is normal or if $J(\nu)$ is tangent. Precisely, we have this results which has to be compared with \cite[Theorem 5.2]{HC}.
\begin{prop}\label{prop3}
Let $(M^n,g)$ be a hypersurface of a Riemannian metallic manifold $(\widetilde{M},\widetilde{g},J)$ with unit normal vector field $\nu$. Let $P:TM\longrightarrow TM$, $V\in\Gamma(TM)$ and $f\in\mathcal{C}^{\infty}(M)$ induced by $J$. If $V$ has no zero and $J(V)$ is normal, or if $J(\nu)$ is tangent, then $M$ is minimal if and only if 
$$\langle \div(P),V \rangle=0.$$
\end{prop}
{\it Proof:} First, assume that $J(V)$ is normal. From Proposition \ref{prop1}, for any $X\in\Gamma(TM)$, we have $A(X)=-\frac{P(\nabla_XV)}{\|V\|^2}=-\frac{P(\nabla_XV)}{q}$. Let $\{e_1,\cdots,e_n\}$ be a local orthonormal frame of $TM$. We have
\beQ
nH&=&\sum_{i=1}^n\langle A(e_i),e_i\rangle\\
&=&-\frac{1}{q}\sum_{i=1}^n\langle P(\nabla_{e_i}V),e_i\rangle\\
&=&-\frac{1}{q}\sum_{i=1}^n\langle \nabla_{e_i}V,P(e_i)\rangle\\
&=&-\frac{1}{q}\sum_{i=1}^n\Big(e_i(\langle V,P(e_i)\rangle)-\langle V,\nabla_{e_i}(P(e_i)) \rangle\Big)\\
&=&-\frac{1}{q}\sum_{i=1}^ne_i(\langle P(V),e_i\rangle)+\langle V,\div(P) \rangle\\
&=&\frac{1}{q}\sum_{i=1}^n\langle V,\div(P) \rangle,
\eeQ
where we have used that $P$ is symmetric and $P(V)=0$. Hence, $M$ is minimal if and only if $\langle \div(P),V \rangle=0$.\\
If $J(\nu)$ is tangent, the proof is ananlogous with the only difference that  $ A(X)=-\frac{P(\nabla_XV)-p\nabla_XV}{\|V\|^2}$ and $P(V)=pV$. The conclusion is the same. \hfill$\square$
\section{Fundamental theorem of submanifolds in product spaces via metallic structures}\label{secpq}
\subsection{Main result}
Now, we assume that $\widetilde{M}$ is the product space $\mathbb{M}^{n_1}(c_1)\times\mathbb{M}^{n_2}(c_2)$, where $\mathbb{M}^{n_1}(c_1)$ and $\mathbb{M}^{n_2}(c_2)$ are the simply connected real space form of respective dimensions $n_1$ and $n_2$ and curvatures $c_1$ and $c_2$. For any positive integers $p$ and $q$, this space is endowed with a canonical metallic structure inherted form the product structure $F$ defined by
\function{F}{T\widetilde{M}=T\mathbb{M}^{n_1}(c_1)\oplus T\mathbb{M}^{n_2}(c_2)}{T\widetilde{M}}{X_1+ X_2}{X_1-X_2.}
We denote by $\pi_1=\frac{1}{2}(\id+F)$ and $\pi_2=\frac{1}{2}(\id-F)$ the projections, respectively on  $T\mathbb{M}^{n_1}(c_1)$ and $T\mathbb{M}^{n_2}(c_2)$. It is a well known fact that the curvature tensor $\widetilde{R}$ of $\widetilde{M}$ is given by
\begin{eqnarray}\label{extcurvature}
\widetilde{R}(X,Y)Z=\sum_{i=1}^2c_i\bigl[\langle\pi_i Y,\pi_iZ\rangle\pi_iX-\langle\pi_i X,\pi_iZ\rangle\pi_iY\bigr].
\end{eqnarray}
Now, as we have seen before, the projections $\pi_1$ and $\pi_2$ can be expressed in term of the metallic structure $J$. Namely, we have
\beqt\label{linkpiJ}
\pi_1=\frac{\sigma_{p,q}}{2\sigma_{p,q}-p}\id-\frac{1}{2\sigma_{p,q}-p}J\quad\text{and}\quad \pi_2=\frac{\sigma_{p,q}-p}{2\sigma_{p,q}-p}\id+\frac{1}{2\sigma_{p,q}-p}J.
\eeqt

Hence, if $M$ is a submanifold of $\widetilde{M}$ then $J$ induces the existence of the four operators $P,Q,R,S$ defined above and the the Gauss, Codazzi and Ricci equations are 

\beq\label{Gauss}
R(X,Y)Z&=&\frac{c_1}{(2\sigma_{p,q}-p)^2}\bigg( \sigma_{p,q}^2(\langle Y,Z\rangle X-\langle X,Z\rangle Y) -\sigma_{p,q}(\langle Y,Z\rangle PX-\langle X,Z\rangle PY)\bigg)\nonumber\\
&&+\frac{c_1}{(2\sigma_{p,q}-p)^2}\bigg(\langle PY,Z\rangle PX-\langle PX,Z\rangle PY-\sigma_{p,q}(\langle PY,Z\rangle X-\langle PX,Z\rangle Y)\bigg)\nonumber\\
&&+\frac{c_2}{(2\sigma_{p,q}-p)^2}\bigg( (\sigma_{p,q}-p)^2(\langle Y,Z\rangle X-\langle X,Z\rangle Y) +(\sigma_{p,q}-p)(\langle Y,Z\rangle PX-\langle X,Z\rangle PY)\bigg)\nonumber\\
&&+\frac{c_2}{(2\sigma_{p,q}-p)^2}\bigg(\langle PY,Z\rangle PX-\langle PX,Z\rangle PY+(\sigma_{p,q}-p)(\langle PY,Z\rangle X-\langle PX,Z\rangle Y)\bigg)\nonumber\\
&&+A_{B(Y,Z)}X-A_{B(X,Z)}Y,
\eeq

\beq\label{Codazzi}
(\nabla_XB)(Y,Z)-(\nabla_YB)(X,Z)&=&\frac{c_1}{(2\sigma_{p,q}-p)^2}\bigg( \langle PY,Z\rangle QX-\langle PX,Z\rangle QY-\sigma_{p,q}(\langle Y,Z\rangle QX-\langle X,Z\rangle QY)\bigg)\nonumber\\
&&+\frac{c_2}{(2\sigma_{p,q}-p)^2}\bigg( \langle PY,Z\rangle QX-\langle PX,Z\rangle QY+(\sigma_{p,q}-p)(\langle Y,Z\rangle QX-\langle X,Z\rangle QY)\bigg),
\eeq

\beq\label{Ricci}
R^{\perp}(X,Y)\nu&=&\frac{c_1+c_2}{(2\sigma_{p,q}-p)^2}\bigg(\langle QY,\nu\rangle QX-\langle QX,\nu\rangle QY\bigg)+B(A_{\nu}Y,Z)-B(A_{\nu}X,Z).
\eeq

Now, we consider $(M^n,g)$ a Riemannian manifold and we also consider a vector bundle $E$ of rank $m$ over $M$ endowed with a metric $g_E$ and a compatible connection $\nabla^E$. We will denote respecively by $\nabla$, $R$ and $R^E$ the Levi-Civita connection of $g$. Moreover let $B:TM\times TM\longrightarrow E$ be a $(2,1)$-symmetric tensor and $P: TM \longrightarrow TM$, $Q: TM \longrightarrow E$, $R:E \longrightarrow TM$ and $S:E \longrightarrow E$ four operators. We give the following definition
\begin{definition}\label{compeq}
We say that $(M,g,E,g^E,\nabla^E,B,P,Q,R,S)$ satisfies the compatiblity equations associated with $(\mathbb{M}^{n_1}(c_1)\times\mathbb{M}^{n_2}(c_2),J)$ if 
\begin{enumerate}
\item $n+m=n_1+n_2$.
\item The map $P$ is symmetric with respect to $g$ and $S$ is symmetric with respect to $g^E$.
\item The maps $Q$ and $R$ are dual with respect to the metric $g\oplus g^E$ over $TM\oplus E$, that is for any $\in TM$ and any $\nu\in E$.
$$g(X,R\nu)=g^E(QX,\nu),$$
\item The rank of the maps $\dfrac{\sigma_{p,q}}{2\sigma_{p,q}-p}\id -\dfrac{1}{2\sigma_{p,q}-p}(P+Q+R+S)$ and\\ $\dfrac{\sigma_{p,q}-p}{2\sigma_{p,q}-p}\id +\dfrac{1}{2\sigma_{p,q}-p}(P+Q+R+S)$ are $n_1$ and $n_2$ respectively.
\item Equations \eqref{relation1.1}-\eqref{relation1.4} are satisfied,
\item Equations \eqref{relation2.1}-\eqref{relation2.4} are satisfied
\item The Gauss, Codazzi and Ricci equations \eqref{Gauss}, \eqref{Codazzi} , \eqref{Ricci} are satisfied.
\end{enumerate}
\end{definition}
As we will see in the following result, these compatibility equations are exactly the necessary and sufficient condition to be immersed isometrically into $\mathbb{M}^{n_1}(c_1)\times\mathbb{M}^{n_2}(c_2)$ with given nornal bundle, second fundamental form and metallic structure over $M$. Namely, we have the following.
\begin{thm}\label{thm1}
Let $(M^n,g)$ be a simply connected Riemannian manifold and $E$ be a $m$-dimensional vector bundle over $M$ endowed with a metric $g^E$ and a compatible conection $\nabla^E$. Moreover, let $B:TM\times TM\longrightarrow E$ be a $(2,1)$-symmetric tensor and $P: TM \longrightarrow TM$, $Q: TM \longrightarrow E$, $R:E \longrightarrow TM$ and $S:E \longrightarrow E$ are four operators. If $(M,g,E,g^E,\nabla^E,B,P,Q,R,S)$ satisfies the compatiblity equations associated with $(\mathbb{M}^{n_1}(c_1)\times\mathbb{M}^{n_2}(c_2),J)$ then, there exists an isometric immersion $\varphi:M\longrightarrow \mathbb{M}^{n_1}(c_1)\times\mathbb{M}^{n_2}(c_2)$ such that the normal bundle of $M$ for this immersion is isometric to $E$ and so that the second fundamental form $II$ and the normal connexion $\nabla^{\perp}$ are given by $B$ and $\nabla^E$. Precisely, there exists a vector bundle isometry $\widetilde{\varphi}: E\longrightarrow T^{\perp}\varphi(M)$ so that
$$II=\widetilde{\varphi}\circ B,$$
$$\nabla^{\perp}\widetilde{\varphi}=\widetilde{\varphi}\nabla^E.$$
Moreover, we have
$$J(\varphi_*X)=\varphi_*(PX)+\widetilde{\varphi}(QX),$$
$$J(\widetilde{\varphi}\nu)=\varphi_*(R\nu)+\widetilde{\varphi}(S\nu),$$
and this immersion is unique up to an isometry of $\mathbb{M}^{n_1}(c_1)\times\mathbb{M}^{n_2}(c_2)$.
\end{thm}
{\it Proof:} We define the following eight operators
$$\left\{
\begin{array}{l}
f_1=\dfrac{\sigma_{p,q}}{2\sigma_{p,q}-p}\id_{TM}-\dfrac{1}{2\sigma_{p,q}-p}P\\ 
f_2=\dfrac{\sigma_{p,q}-p}{2\sigma_{p,q}-p}\id_{TM}+\dfrac{1}{2\sigma_{p,q}-p}P\\ 
h_1=-h_2=-\dfrac{1}{2\sigma_{p,q}-p}Q \\ 
s_1=-s_2= -\dfrac{1}{2\sigma_{p,q}-p}R\\ 
t_1= \dfrac{\sigma_{p,q}}{2\sigma_{p,q}-p}\id_{E}-\dfrac{1}{2\sigma_{p,q}-p}S\\ 
t_2= \dfrac{\sigma_{p,q}-p}{2\sigma_{p,q}-p}\id_{E}+\dfrac{1}{2\sigma_{p,q}-p}S.
\end{array}
\right.$$
It is clear from the definition of these eigth operators that 
$$f_1,f_2:TM\lgra TM,\ h_1,h_2: TM\lgra E,\ s_1,s_2:E\lgra TM,\ t_1,t_2:E\lgra E$$ and that
$$f_1+f_2=\id_{TM},\ h_1+h_2=0,\ s_1+s_2=0,\ t_1+t_2=\id_{E}.$$
Moreover, we have the following lemma.
\begin{lemma}
The operators $f_1,f_2,h_1,h_2,s_1,s_2,t_1,t_2$ satisfy the following relations:\\

\begin{align}
&f_i\circ f_j+s_i\circ h_j=\delta_i^jf_i,&\nonumber\\
&t_i\circ t_j+h_i\circ s_j=\delta_i^jt_i,&\nonumber \\
&f_i\circ s_j+s_i\circ t_j=\delta_i^js_i,&\nonumber\\
&h_i\circ f_j+t_i\circ h_j=\delta_i^jh_i,&\nonumber\\
&\nabla_X(f_iY)-f_i(\nabla_XY)=A_{h_iY}X+s_i(B(X,Y)),& \nonumber\\
&\overline{\nabla}_X(h_iY)-h_i(\nabla_XY)=t_i(B(X,Y))-B(X,f_iY),& \nonumber\\
&\overline{\nabla}_X(t_i\nu)-t_i(\overline{\nabla}_X\nu)=-B(s_i\nu,X)-h_i(A_{\nu}X),&\nonumber\\
&\nabla_X(s_i\nu)-s_i(\nabla^{\perp}_X\nu)=-f_i(A_{\xi}X)+A_{t_i\nu}X,&\nonumber
\end{align}
for $i,j\in\{1,2\}$, any $X,Y\in TM$ and $\nu\in E$.
\end{lemma}
{\it Proof:} The proof of this lemma is elementary. We will give only the proof of two relation to give an idea of the computations (which are straightforward). First, we have
\beQ
f_1\circ f_2+s_1\circ h_2&=&\Big(\dfrac{\sigma_{p,q}}{2\sigma_{p,q}-p}\id_{TM}-\dfrac{1}{2\sigma_{p,q}-p}P\big)\circ \big( \dfrac{\sigma_{p,q}-p}{2\sigma_{p,q}-p}\id_{TM}+\dfrac{1}{2\sigma_{p,q}-p}P\Big)\\
&&+\Big( -\dfrac{1}{2\sigma_{p,q}-p}R\Big)\circ \Big( \dfrac{1}{2\sigma_{p,q}-p}Q\Big)\\
&=&\dfrac{\sigma_{p,q}(\sigma_{p,q}-p)}{(2\sigma_{p,q}-p)^2}\id_{TM}+\frac{p}{(2\sigma_{p,q}-p)^2}P-\frac{1}{(2\sigma_{p,q}-p)^2}P^2+\frac{1}{(2\sigma_{p,q}-p)^2}R\circ Q\\
&=&\frac{1}{(2\sigma_{p,q}-p)^2}(-P^2+pP-R\circ Q-(\sigma_{p,q}-p)\sigma_{p,q}\id_{TM})\\
&=&0,
\eeQ
where we have used that $(\sigma_{p,q}-p)\sigma_{p,q}=q$ and \eqref{relation1.1}, i.e., $P^2+R\circ Q=pP+q\id_{TM}$. The proof of the other identites of the same type are analogous and use \eqref{relation1.1} to \eqref{relation1.4}.\\
Moreover, we have
\beQ
\nabla_X(f_1Y)-f_1(\nabla_XY)&=&\frac{\sigma_{p,q}}{(2\sigma_{p,q}-p)}\nabla_XY-\frac{1}{(2\sigma_{p,q}-p)}\nabla_X(PY)-\frac{\sigma_{p,q}}{(2\sigma_{p,q}-p)}\nabla_XY+\frac{1}{(2\sigma_{p,q}-p)}P(\nabla_XY)\\
&=&-\frac{1}{(2\sigma_{p,q}-p)}\Big(\nabla_X(PY)-P(\nabla_XY)\Big)\\
&=&-\frac{1}{(2\sigma_{p,q}-p)}\Big(  A_{QY}X-R(B(X,Y))\Big)\\
&=&A_{h_1Y}X-s_1(B(X,Y)).
\eeQ
The other identites are proven in the same way, using relations \eqref{relation2.1} to \eqref{relation2.4}. \hfill$\square$\\ \\
Finally, we can see easily that the Gauss, Codazzi and Ricci Equations \eqref{Gauss}, \eqref{Codazzi} and \eqref{Ricci} can be written as follows with the definition of $f_1,f_2,h_1,h_2,s_1,s_2,t_1,t_2$ form $P,Q,R,S$:
\begin{align}
&R(X,Y)Z=\sum_{i=1}^mc_i\bigg[ \left\langle f_iY,Z\right\rangle f_iX-\left\langle f_iX,Z\right\rangle f_iY \bigg]+A_{B(Y,Z)}X-A_{B(X,Z)}Y,&\nonumber\\
&(\nabla_XB)(Y,Z)-(\nabla_YB)(X,Z)=\sum_{i=1}^mc_i\bigg[ \left\langle f_iY,Z\right\rangle h_iX-\left\langle f_iX,Z \right\rangle h_iY\bigg],&\nonumber\\
&\overline{R}(X,Y)\nu=\sum_{i=1}^mc_i\bigg[ \left\langle h_iY,\nu\right\rangle h_iX-\left\langle h_iX,\nu\right\rangle h_iY\bigg]+B(A_{\nu}Y,X)-B(A_{\nu}X,Y).\nonumber&
\end{align}
We conclude by using Theorem 3.2 of \cite{LR} to conclude. Indeed, we have shown that $(M,g,E,g^E,\nabla^E,B,f_1,f_2,h_1,h_2,s_1,s_2,t_1,t_2)$ satisfies the compatibility equations for the isometric immersions into the products space $\mathbb{M}^{n_1}(c_1)\times\mathbb{M}^{n_2}(c_2)$ defined in \cite{LR}. Hence, by Theorem 3.2 of \cite{LR}, there exists an isometric immersion $\varphi:M\lgra\mathbb{M}^{n_1}(c_1)\times\mathbb{M}^{n_2}(c_2)$ such that the normal bundle of $M$ for this immersion is isomorphic to $E$ and such that the second fundamental form $II$ and the normal connection $\nabla^{\perp}$  are given by $B$ and $\nabla^E$, that is, there exists a vector bundle isometry $\widetilde{\varphi}:E\lgra T^{\perp}\varphi(M)$ so that
$$II=\widetilde{\varphi}\circ B,$$
$$\nabla^{\perp}\widetilde{\varphi}=\widetilde{\varphi}\nabla^E.$$
Moreover, we have for $i=1,2$, 
\beqt\label{phi1}
\pi_i(\varphi_*X)=\varphi_*(f_iX)+\widetilde{\varphi}(h_iX),
\eeqt
\beqt\label{phi2}
\pi_i(\widetilde{\varphi}\nu)=\varphi_*(s_i\nu)+\widetilde{\varphi}(t_i\nu),
\eeqt
where $\pi$ is the projection on $T\mathbb{M}^{n_i}(c_i)$,
and this isometric immersion is unique up to an isometry of $\mathbb{M}^{n_1}(c_1)\times\mathbb{M}^{n_2}(c_2)$. Finally, from the definition of $f_1,f_2,h_1,h_2,s_1,s_2,t_1,t_2$, we get
$$\left\{
\begin{array}{l}
P=\sigma_{p,q}f_2-(\sigma_{p,q}-p)f_1,\\
Q=-2(\sigma_{p,q}-p)h_1=2(\sigma_{p,q}-p)h_2,\\
R=-2(\sigma_{p,q}-p)s_1=2(\sigma_{p,q}-p)s_2,\\
S=\sigma_{p,q}t_2-(\sigma_{p,q}-p)t_1.
\end{array}
\right.$$
Using these relations into \eqref{phi1} and \eqref{phi2}, together with \eqref{linkpiJ} we get finally
$$J(\varphi_*X)=\varphi_*(PX)+\widetilde{\varphi}(QX),$$
$$J(\widetilde{\varphi}\nu)=\varphi_*(R\nu)+\widetilde{\varphi}(S\nu).$$
This concludes the proof of the theorem.\hfill$\square$\\ \\
Now, we can obtain an analogue to the main result of \cite{Ro3} in terms of $(p,q)$-metallic structures. In \cite{Ro3}, we obtain a spinorial characterization of surfaces into the $4$-dimensional products $\mathbb{M}^2(c)\times \RR^2$ and $\mathbb{M}^3(c)\times \RR$ generalizing the results for space forms proven in \cite{BLR}. For all the basics about the spinorial geometry of surfaces into $4$-dimensional spaces, the reader can refer to \cite{BLR,Ro3,NR}. We prove the following.
\begin{thm}\label{corM2R2}
Let $c\in\RR$, $c\neq0$ and $\alpha\in\CC$ such that $4\alpha^2=c$. Let $(M^2,g)$ be an oriented Riemannian surface and $E$ an oriented vector bundle of rank $2$ over $M$ with scalar product $\langle\cdot,\cdot\rangle_E$ and compatible connection $\nabla^E$. We denote by $\Sigma=\Sigma M\otimes\Sigma E$ the twisted spinor bundle. Let $B:TM\times TM\lgra E$ a bilinear symmetric map and 
$$P:TM\lgra TM,\ Q:TM\lgra E,\ R:E\lgra TM\ \text{and}\ S:E\lgra E$$ satisfying equations \eqref{relation1.1}-\eqref{relation2.4}. Moreover we assume that the rank of the maps $\dfrac{\sigma_{p,q}}{2\sigma_{p,q}-p}\id -\dfrac{1}{2\sigma_{p,q}-p}(P+Q+R+S)$ and $\dfrac{\sigma_{p,q}-p}{2\sigma_{p,q}-p}\id +\dfrac{1}{2\sigma_{p,q}-p}(P+Q+R+S)$ are 2 and 2 ({\it resp.} 3 and 1). Then, the two following statements are equivalent
\begin{enumerate}
\item There exists a spinor field $\varphi$ in $\Sigma$  satisfying for all $X\in\Chi(M)$
$$\nabla_X\varphi=\left(\dfrac{\alpha\sigma_{p,q}}{2\sigma_{p,q}-p}X-\dfrac{\alpha}{2\sigma_{p,q}-p}(PX+QX)\right)\cdot\varphi-\frac{1}{2}\zeta(X)\cdot\varphi,$$ such that $\varphi^+$ and $\varphi^-$ vanish nowhere and where $\zeta$ is given by
$$
\zeta(X)=\sum_{j=1}^2e_j\cdot B(e_j,X).
$$
\item There exists a local isometric immersion of $(M^2,g)$ into $P=\M^2(c)\times\RR^2$ ({\it resp.} $\M^3(c)\times\RR$) with $E$ as normal bundle and second fundamental form $B$ such that over $M$ the canonical $(p,q)$-metallic structure is given by $P,Q,R$ and $S$ in the sense of Theorem \ref{thm1}.

\end{enumerate}
\end{thm}
{\it Proof:} This theorem is a direct application of \cite[Theorem 3.1]{Ro3}. Indeed, if there exists a isometric immersion of $(M^2,g)$ into $\M^2(c)\times\RR^2$ ({\it resp.} $\M^3(c)\times\RR$), then, as proved in \cite{Ro3}, there exists a spinor field satisfying 
\begin{equation}\label{eqspin}
\nabla_X\varphi=\frac{\alpha}{2}(X+fX+hX)\cdot\varphi+\zeta(X)\cdot\varphi,
\end{equation}
where $f,h,s,t$ are the four operators induce by the structure product $F$. Hence, with the link between products and $(p,q)$-metallic strctures
$$\frac{2}{2\sigma_{p,q}-p}J-\frac{p}{2\sigma_{p,q}-p}\id,$$
we get 
$$\nabla_X\varphi=\left(\dfrac{\alpha\sigma_{p,q}}{2\sigma_{p,q}-p}X-\dfrac{\alpha}{2\sigma_{p,q}-p}(PX+QX)\right)\cdot\varphi+\zeta(X)\cdot\varphi.$$
Conversely, if all the assumption of point $(2)$ are satisfied, by setting
$$\left\{
\begin{array}{l}
f=\dfrac{2}{2\sigma_{p,q}-p}P-\dfrac{p}{2\sigma_{p,q}-p}\id_{TM}\\
h=\dfrac{2}{2\sigma_{p,q}-p}Q\\
s=\dfrac{2}{2\sigma_{p,q}-p}R\\
t=\dfrac{2}{2\sigma_{p,q}-p}S-\dfrac{p}{2\sigma_{p,q}-p}\id_{E},
\end{array}
\right.$$
 straigthforward computations show that the assumptions of Theorem 3.1 of \cite{Ro3} are satisfied and so $M$ is locally isometrically immersed into the {\it ad hoc} product space with $E$ as normal bundle, $B$ as second fundamental form and over $M$, the product structure is given by $f,h,s$ and $t$. But, form the definition of $f,h,s$ and $t$ from $P,Q,R$ and $S$, then $P,Q,R$ and $S$ are necessarily the restriction of the $(p,q)$-metallic structure $J=\frac{p}{2}\id+\frac{2\sigma_{p,q}-p}{2}F$. This concludes the proof.\hfill$\square$
\subsection{Associated families}
In this section, we use Theorem \ref{thm1} to recover the theorem of existence of  associate families of minimal surfaces into the multiproduct $\mathbb{M}^{n_1}(c_1)\times\mathbb{M}^{n_2}(c_2)$ expressed in terms of metallic structure.\\
Let $(\Sigma,g)$ be an oriented Riemannian surface. We denote by $\mathcal{J}$ its complex structure, that is, the rotation of angle $\frac{\pi}{2}$ on $TM$. For any $\theta\in\RR$, we set $\mathcal{R}_{\theta}=\cos(\theta)I+\sin(\theta)\mathcal{J}$. First, we have the following proposition. 
\begin{prop}\label{prop1}
Assume that $(\Sigma,g,E,g^E,\nabla^E,B,P,Q,R,S)$ satisfies the compatibility equation for $(\mathbb{M}^{n_1}(c_1)\times\mathbb{M}^{n_2}(c_2),J)$ and that $B$ is trace-free for any $\nu\in E$, then $(\Sigma,g,E,g^E,\nabla^E,B_{\theta},f_{i,\theta},h_{i,\theta},t_{i,\theta})$ also satisfies the compatibility equations for $(\mathbb{M}^{n_1}(c_1)\times\mathbb{M}^{n_2}(c_2),J)$, where
\begin{align}
&B_{\theta}(X,Y)=B(X,\mathcal{R}_{\theta}^{-1}Y),&\nonumber\\
&P_{\theta}=\mathcal{R}_{\theta}\circ P\circ\mathcal{R}^{-1}_{\theta},&\nonumber\\
&Q_{\theta}=Q\circ\mathcal{R}^{-1}_{\theta},&\nonumber\\
&R_{\theta}=\mathcal{R}_{\theta}\circ R,&\nonumber\\
&S_{\theta}=S.&\nonumber
\end{align}
Moreover, $B_{\theta}$ is also trace-free for any $\nu\in E$.
\end{prop}
{\it Proof:} From $P,Q,R$ and $S$, we have defined $f_1,f_2,h_1,h_2,s_1,s_2,t_1,t_2$. Moreover, from $P_{\theta},Q_{\theta},R_{\theta}$ and $S_{\theta}$, we define in the same way $f_{1,\theta},f_{2,\theta},h_{1,\theta},h_{2,\theta},s_{1,\theta},s_{2,\theta},t_{\theta}1,t_{2,\theta}$. It is clear from the definition of all these operators that we have for $i=1,2$
\begin{align}
&f_{i,\theta}=\mathcal{R}_{\theta}\circ f_i\circ\mathcal{R}^{-1}_{\theta},&\nonumber\\
&h_{i,\theta}=h_i\circ\mathcal{R}^{-1}_{\theta},&\nonumber\\
&s_{i,\theta}=\mathcal{R}_{\theta}\circ s_i,&\nonumber\\
&t_{i,\theta}=t_i.&\nonumber
\end{align}
Hence, as we have seen in the proof of Theorem \ref{thm1}, $(M,g,E,g^E,\nabla^E,B,f_1,f_2,h_1,h_2,s_1,s_2,t_1,t_2)$ satisfies the compatibility equations for the isometric immersions into the products space $\mathbb{M}^{n_1}(c_1)\times\mathbb{M}^{n_2}(c_2)$ defined in \cite{LR} and using Proposition 4.1 in \cite{LR}, we get that $B_{\theta}$ is trace-free and $(M,g,E,g^E,\nabla^E,B,f_{1,\theta},f_{2,\theta},h_{1,\theta},h_{2,\theta},s_{1,\theta},s_{2,\theta},t_{\theta}1,t_{2,\theta})$ also satisfied the same compatibility conditions. Finally, since we have
$$\left\{
\begin{array}{l}
P_{\theta}=\sigma_{p,q}f_{2,\theta}-(\sigma_{p,q}-p)f_{1,\theta},\\
Q_{\theta}=-2(\sigma_{p,q}-p)h_{1,\theta}=2(\sigma_{p,q}-p)h_{2,\theta},\\
R_{\theta}=-2(\sigma_{p,q}-p)s_{1,\theta}=2(\sigma_{p,q}-p)s_{2,\theta},\\
S_{\theta}=\sigma_{p,q}t_{2,\theta}-(\sigma_{p,q}-p)t_{1,\theta},
\end{array}
\right.$$
we get immediately that $(\Sigma,g,E,g^E,\nabla^E,B_{\theta},f_{i,\theta},h_{i,\theta},t_{i,\theta})$ satisfies the compatibility conditions of Definition \ref{compeq}. This concludes the proof. \hfill$\square$\\ \\
Now, using this proposition and Theorem \ref{thm1}, we can deduce easily the following result.
\begin{thm}\label{thm2}
Let $\Sigma$ be a simply connected surface and $\varphi:M\lgra \mathbb{M}^{n_1}(c_1)\times\mathbb{M}^{n_2}(c_2)$ be a minimal isometric immersion with normal bundle $E$, second fundamental form $B$ and normal connection $\nabla^{\perp}$. Let $P$, $Q$, $R$ and $S$ be the $(1,1)$-tensors induced by the metallic structure $J$. Let $p_0\in\Sigma$. Then, there exists a unique family $(\varphi_{\theta})_{\theta\in\RR}$ of minimal isometric immersions $\varphi_{\theta}:\Sigma\lgra\mathbb{M}^{n_1}(c_1)\times\mathbb{M}^{n_2}(c_2)$ so that
\begin{enumerate}[(i)]
\item $\varphi_{\theta}(p_0)=\varphi(p_0)$ and $d({\varphi_{\theta}})_{p_0}=(d\varphi)_{p_0}$,
\item the metric induced by $\varphi$ and $\varphi_{\theta}$ are the same,
\item the second fundamental form of $\varphi_{\theta}(\Sigma)$ in $\mathbb{M}^{n_1}(c_1)\times\mathbb{M}^{n_2}(c_2)$ is given by $B_{\theta}(X,Y)=B(R_{\theta}X,Y)$, for any $X,Y\in\Gamma(T\Sigma)$,
\item for any $X\in\Gamma(T\Sigma)$ and $\xi\in\Gamma(E)$,
$$J(d\varphi_{\theta}X)=d\varphi_{\theta}(P_{\theta}X)+Q_{\theta}X\quad \text{and}\quad J(\xi)=d\varphi_{\theta}(R_{\theta}\xi)+S_{\theta}\xi.$$
\end{enumerate}
Moreover, $\varphi_0=\varphi$ and the family $(\varphi_{\theta})_{\theta\in\RR}$ is continuous with respect to $\theta$.
\end{thm}
{\it Proof:} Since $\Sigma$ is a minimal surface in $\mathbb{M}^{n_1}(c_1)\times\mathbb{M}^{n_2}(c_2)$, then $(\Sigma,g,E,g^E,\nabla^E,B,P,Q,R,S)$ satisfy the compatibility equations and so, by Proposition \ref{prop1}, $(\Sigma,g,E,g^E,\nabla^E,B_{\theta},P_{\theta},Q_{\theta},R_{\theta},S_{\theta})$ also satisfies the compatibility equations. Using, Theorem \ref{thm1}, we know that there exists an isometric immersion from $\Sigma$ into $\mathbb{M}^{n_1}(c_1)\times\mathbb{M}^{n_2}(c_2)$ which satisfy the point $(iii)$ and $(iv)$. Moreover, this immersion is also minimal since $B_{\theta}$ is trace-free by Proposition \ref{prop1}. Finally the point $(i)$ is clear by construction, the point $(ii)$ also because both induced metric are $g$ and the continuity of the family is ensured by the construction of the immersions in \cite{LR}.
\subsection{Examples}
We finish this note with some examples. We consider the product of two spheres $M=\SSS^{n_1}(c_1)\times\SSS^{n_2}(c_2)$ endowed the product metric. We denote by $r_i=\frac{1}{\sqrt{c_1}}$ the radius of the sphere of curvature $c_1$. Obviously, $M$ can be canonically isometrically embedded into the Euclidean space $\RR^{n_1+n_2+2}=\RR^{n_1+1}\times\RR^{n_2+1}$.  For more convenience and compactness, we will denote $(x_1,\cdots,x_{n_1+1},y_1,\cdots,y_{n_2+1})$ by $(x,y)$ with $x=(x_1,\cdots,x_{n_1+1})$ and $y=(y_1,\cdots,y_{n_2+1})$. Hence $\SSS^{n_1}(c_1)\times\SSS^{n_2}(c_2)$ is defined by
$$\SSS^{n_1}(c_1)\times\SSS^{n_2}(c_2)=\left\{(x,y)\in \RR^{n_1+n_2+2}\ |\ \sum_{i=1}^{n_1+1}x_i^2=r_1^2\ \text{and}\ \sum_{i=1}^{n_2+1}y_i^2=r_2^2\right\}.$$
We denote by $E$ the normal bundle of this immersion, $g^E$ the induced normal metric and $\nabla^E$ the induced normal connection. Note that the normal connection is flat. The canonical embedding is just the inclusion map. Moreover, a vector $Z$ tangent to $M$ at the point $(x,y)$ is of the form $Z=(X,Y)$ with $\displaystyle\sum_{i=1}^{n_1+1}x_iX_i=0$ and $\displaystyle\sum_{i=1}^{n_2+1}y_iY_i=0$. The normal bundle has a global orthonormal frame $\{N_1,N_2\}$ with $N_1=\frac{1}{r_1}(x,0)$ and $N_2=\frac{1}{r_2}(0,y)$.

First, we consider the following $(p,q)$-metallic structure on $\RR^{n_1+n_2+2}=\RR^{n_1+1}\times\RR^{n_2+1}$:
\function{\widetilde{P}}{\RR^{n_1+n_2+2}=\RR^{n_1+1}\times\RR^{n_2+1}}{\RR^{n_1+n_2+2}=\RR^{n_1+1}\times\RR^{n_2+1}}{(x_1,\cdots,x_{n_1+1},y_1,\cdots,y_{n_2+1})}{(\sigma_{p,q}x_1,\cdots,\sigma_{p,q}x_{n_1+1},(p-\sigma_{p,q})y_1,\cdots,(p-\sigma_{p,q})y_{n_2+1})}.
Clearly, $\widetilde{P}$ is a $(p,q)$-metallic structure (see \cite{HC} for the details) and as we have seen in Section \ref{secsubmdfpq}, $\widetilde{P}$ induces the existence of the four operators $P,Q,R$ and $S$. It is obvious to see that $Q=0$, $R=0$ and that $P$ and $R$ are given by:
$$P(Z)=P(X,Z)=(\sigma_{p,q}X,(p-\sigma_{p,q})Y)$$
and 
$$S(N_1)=\sigma_{p,q}N_1,\quad S(N_2)=(p-\sigma_{p,q})N_2.$$
Note that $P$ and $S$ are metallic structures on $TM$ and $E$, respectively. Moreover, the curavture tensor of $\SSS^{n_1}(c_1)\times\SSS^{n_2}(c_2)$ is given by
\beQ
R(X,Y)Z&=&\frac{c_1}{(2\sigma_{p,q}-p)^2}\bigg( \sigma_{p,q}^2(\langle Y,Z\rangle X-\langle X,Z\rangle Y) -\sigma_{p,q}(\langle Y,Z\rangle PX-\langle X,Z\rangle PY)\bigg)\nonumber\\
&&+\frac{c_1}{(2\sigma_{p,q}-p)^2}\bigg(\langle PY,Z\rangle PX-\langle PX,Z\rangle PY-\sigma_{p,q}(\langle PY,Z\rangle X-\langle PX,Z\rangle Y)\bigg)\nonumber\\
&&+\frac{c_2}{(2\sigma_{p,q}-p)^2}\bigg( (\sigma_{p,q}-p)^2(\langle Y,Z\rangle X-\langle X,Z\rangle Y) +(\sigma_{p,q}-p)(\langle Y,Z\rangle PX-\langle X,Z\rangle PY)\bigg)\nonumber\\
&&+\frac{c_2}{(2\sigma_{p,q}-p)^2}\bigg(\langle PY,Z\rangle PX-\langle PX,Z\rangle PY+(\sigma_{p,q}-p)(\langle PY,Z\rangle X-\langle PX,Z\rangle Y)\bigg).\nonumber
\eeQ
\subsubsection{Isometric immersion into $\SSS^{n_1+1}(c_1)\times\SSS^{n_2+1}(c_2)$.}
Now, we forget the immersion of $\SSS^{n_1}(c_1)\times\SSS^{n_2}(c_2)$ but only consider $(E,g^E,\nabla^E)$ as a rank $2$ vector bundle over $\SSS^{n_1}(c_1)\times\SSS^{n_2}(c_2)$ with given metric and compatible connection. We consider the same operators $P,Q,R$ and $S$ and we set $B=0$. By straight forward computations, we see that all the relations of Proposition \ref{allrelations} are satisfied. Moreover, since $B=0$, $Q=0$ and the curvature associated with $\nabla^E$ is zero, then the Codazzi and Ricci equations are trivial. Finally, since $B=0$, the expression of the curvature tensor of $\SSS^{n_1}(c_1)\times\SSS^{n_2}(c_2)$ gives immediately the Gauss equations. Finally, we can apply Theorem \ref{thm1} and recover the isometric immersion of $\SSS^{n_1}(c_1)\times\SSS^{n_2}(c_2)$ into $\SSS^{n_1+1}(c_1)\times\SSS^{n_2+1}(c_2)$ with $E$ as normal bundle, vanishing second fundamental form and so that $P$ and $S$ are the restrictions over $TM$ and $E$ respectively of the canonical $(p,q)$-metallic structure of $\SSS^{n_1+1}(c_1)\times\SSS^{n_2+1}(c_2)$.
\subsubsection{Isometric immersion into $\SSS^{n_1}(c_1)\times\RR^{n_2+1}$.}
Now, we consider the same operator $P$ on $\SSS^{n_1}(c_1)\times\SSS^{n_2}(c_2)$ and we define $V=0$ and $f=p-\sigma_{p,q}$. Moreover, we define $A:TM\lgra TM$ by
$$AZ=\sqrt{c_2}\left(\frac{\sigma_{p,q}}{2\sigma_{p,q}-p}Z-\frac{1}{2\sigma_{p,q}-p}PZ\right).$$
We have the following elementary lemma.
\begin{lemma}
All the relations of Proposition \ref{allrelationshyp} as well as the Gauss and Codazzi euqations are satisfied.
\end{lemma}
{\it Proof:} We have already say that $P$ is a $(p,q)$-metallic structure, so relation \eqref{relation1.5hyp} is satisfied. Moreover, since $V=0$, relation \eqref{relation1.1hyp} and \eqref{relation2.1hyp} also satisfied as well as \eqref{relation1.2hyp} and \eqref{relation2.4hyp} which are trivial. Moreover,  $f=p-\sigma_{p,q}$, then $f$ satisfies $f^2=pf=q$ which is relation \eqref{relation1.4hyp} since $V=0$. Moreover, we have
\beQ
fAZ-P(AZ)&=&f\sqrt{c_2}\left(\frac{\sigma_{p,q}}{2\sigma_{p,q}-p}Z-\frac{1}{2\sigma_{p,q}-p}PZ\right)-\sqrt{c_2}P\left(\frac{\sigma_{p,q}}{2\sigma_{p,q}-p}Z-\frac{1}{2\sigma_{p,q}-p}PZ\right)\\
&=&\frac{\sqrt{c_1}}{2\sigma_{p,q}-p}\left( f\sigma_{p,q}Z-fPZ-\sigma_{p,q}PZ+P^2Z \right)\\
&=&\frac{\sqrt{c_1}}{2\sigma_{p,q}-p}\left( P^2Z -pPZ-qZ\right)\\
&=&0,
\eeQ
and so \eqref{relation2.3hyp} is satisfied.
Moreover, since $P$ is parallel, then $A$ is also parallel and since $V=0$, the Codazzi equation is trivial. Finally, the expression of the curvature of $\SSS^{n_1}(c_1)\times\SSS^{n_2}(c_2)$ give immediately the Gauss equation from the definition of $A$.
\hfill$\square$\\ \\
Thus, the compatibilty equations are satisfied to recover an isometric immersion of $\SSS^{n_1}(c_1)\times\SSS^{n_2}(c_2)$ into  $\SSS^{n_1}(c_1)\times\RR^{n_2+1}$ as a hypersurface with shape operator given by $A$ and so that the restriction of the canonical $(p,q)$-metallic structure of $\SSS^{n_1}(c_1)\times\RR^{n_2+1}$ is given by $P$ and $f$.
\section{Complex metallic structures}
In this section, we will define complex metallic structures, which are in some sense the analogue for complex structures of what metallic structures are for product structures. After giving the definition and basic properties of these new structures, we study submanifolds of Riemannian manifolds admitting complex metallic structures. In particular, we prove in this context similar results that those proved above for $(p,q)$-metallic structures.
\subsection{Definition}
Let $a,b,$ be two positive real number so that $a<2\sqrt{b}$. We consider the second degree equation $x^2+ax+b=0$. This equation has complex conjugated solutions $\dfrac{-a\pm i\sqrt{4b-a^2}}{2}$. In the sequel, we will denote $\delta=\sqrt{4b-a^2}$.
\begin{definition}
\begin{itemize}
\item
A $(a,b)$-complex metallic structure on a manifold $M$ is a $(1,1)$-tensor $J$ on $M$ satisfying the relation 
$$J^2+aJ+b\id_{TM}=0.$$
\item
If $M$ is endowed with a Riemannian metric $g$, then a $(a,b)$-complex metallic structure $J$ on $M$ is said to be Riemannian if $J$ satisfies for any $X,Y\in\Gamma(TM)$,
$$g(JX,Y)=-g(X,JY)-ag(X,Y).$$
\item A Riemannian $(a,b)$-complex metallic structure $J$ on $(M,g)$ is said parallel if $J$ is parallel with respect to the Levi-Civita connection of $g$.
\end{itemize}
\end{definition}
As we have seen in section \ref{secmetallic}, $(p,q)$-structures are in correspondence with products structures. We have a comparable result for $(a,b)$-complex metallic structures and complex structures.
\begin{prop}\label{equivstructurecomp}Let $M$ be a smooth manifold.  Then, we have
\begin{enumerate}
\item Every almost complex structure $\mathcal{J}$ on $M$ induces two $(a,b)$-metallic structures $J_1,J_2$ on $M$ defined by
$$J_+=-a\id+\delta\mathcal{J}\quad\text{and}\quad J_-=-a\id-\delta\mathcal{J}.$$
Moreover, if $g$ is a Riemannian metric on $M$ and if $(M,g,\mathcal{J})$ is K\"ahler, then $J_1$ and $J_2$ are Riemannian parallel $(a,b)$-metallic structures.
\item Every $(a,b)$-metallic structure $J$ on $M$ induces the existence of two almost complex structures $\mathcal{J}_{\pm}$ defined by
$$\mathcal{J}_{\pm}=\pm\left(\frac{2}{\delta}J+\frac{a}{\delta}\id\right).$$
\item If $M$ carries a $(a,b)$-metallic structure, then $M$ is even-dimensional.
\end{enumerate}
\end{prop}
{\it Proof:} $(1)$ Let $\mathcal{J}$ be an almost complex structure on $M$. We set $J_{\pm}=-a\id\pm\delta\mathcal{J}$. Then, we have
\beQ
J_{\pm}^2+aJ_{\pm}+b\id&=&\left(-\frac{a}{2}\id\pm\frac{\delta}{2}\mathcal{J}\right)^2+a\left(-\frac{a}{2}\id\pm\frac{\delta}{2}\mathcal{J}\right)+b\id\\
&=&\frac{a^2}{4}\id\mp\frac{a}{2}\delta\mathcal{J}-\frac{\delta^2}{4}\id-\frac{a^2}{2}\id\pm\frac{a}{2}\delta\mathcal{J}+b\id\\
&=&\left(b-\frac{a^2}{4}-\delta^2\right)\id\\
&=&0,
\eeQ
since $\delta=\sqrt{b-4a^2}$. Moreover, if $\mathcal{J}$ is compatible with $g$, we have
\beQ
g(J_{\pm}X,Y)&=&g\left(-\frac{a}{2}\pm\frac{\delta}{2}\mathcal{J}X,Y\right)\\
&=&-\frac{a}{2}g(X,Y)\mp \frac{\delta}{2}g(X,\mathcal{J}Y)\\
&=&-g\left(X,-\frac{a}{2}\pm\frac{\delta}{2}\mathcal{J}Y\right)-ag(X,Y),
\eeQ
and so $J_{\pm}$ are Riemannian $(a,b)$-complex metallic structures. Finally, if $\mathcal{J}$ is parallel, then $J_{\pm}$ are clearly parallel.\\ \\
$(2)$  Let $J$ be a $(a,b)$-complex metallic structure. We set $\mathcal{J}_{\pm}=\pm\left(\frac{2}{\delta}J+\frac{a}{\delta}\id\right)$. We have
\beQ
\mathcal{J}_{\pm}^2&=&\left(\frac{2}{\delta}J+\frac{a}{\delta}\id\right)^2\\
&=&\frac{4}{\delta^2}J^2+\frac{4a}{\delta^2}J+\frac{a^2}{\delta^2}\id\\
&=&\frac{4}{\delta^2}\left( J^2+aJ+b\id\right)+\frac{1}{\delta^2}(a^2-4b)\id\\
&=&-\id,
\eeQ
where, we have used that $J^2+aJ+b\id=0$ and $\delta^2=4b-a^2$. Hence, $\mathcal{J}_{\pm}$ are almost complex structures. Moreover, if $J$ is a Riemannian $(a,b)$-complex structure, then by definition, we have for any $X,Y\in\Gamma(TM)$. 
$$g(JX,Y)=-g(X,JY)-ag(X,Y).$$
From this, we get
\beQ
g(\mathcal{J}_{\pm}X,Y)&=&\pm g\left(\frac{2}{\delta}JX+\frac{a}{\delta}X,Y\right)\\
&=&\mp\frac{2}{\delta}g(X,JY)\mp\frac{2a}{\delta}g(X,Y)\pm\frac{a}{\delta}g(X,Y)\\
&=&\mp g\left(X,\frac{2}{\delta}JY+\frac{a}{\delta}Y\right)\\
&=&-g(X,\mathcal{J}_{\pm}Y)
\eeQ
and so $\mathcal{\pm}$ is compatible with $g$. Finally, it is clear that if $J$ is parallel, then also $\mathcal{J}_{\pm}$.\\ \\
$(3)$ If $M$ carries a $(a,b)$-complex metallic structure, then $M$ also carries a complex structure by point $(2)$, and so $M$ is necessarily even-dimensional.
\subsection{Submanifolds of complex metallic structures}

Now, let us consider a Riemannian manifold $(M^n,g)$ isometrically immersed into a $(n+m)$-dimensional Riemannian manifold $(\widetilde{M},\widetilde{g})$ endowed with a Riemannin parallel $(a,b)$-complex metallic structure $J$. We denote by $E$ the normal bundle which is equipped with an induced metric $g^E$ and the induced compatible normal connection $\nabla^E$. Then, the complex metallic structure $J$ induces the existence of four operators $P: TM \longrightarrow TM$, $Q: TM \longrightarrow E$, $R:E \longrightarrow TM$ and $S:E \longrightarrow E$, so that with respect to the decomposition $T\widetilde{M}=TM\oplus E$,  $J$ is given over $M$ by
$J=\left( 
\begin{array}{cc}
P&R\\
Q&S
\end{array}\right).$
Then, the operators $P,Q,R$ and $S$ satisfy the following equations:
\begin{prop}\label{allrelationscomp}
For all $X,Y\in TM$ and all $\xi,\nu\in E$, we have
\begin{align}
&P^2+R\circ Q=-aP-b\id_{TM},&\label{relation1.1comp}\\
&Q\circ P+S\circ Q=-aQ,&\label{relation1.2comp}\\
&P\circ R+R\circ S=-aR,&\label{relation1.3comp}\\
&S^2+Q\circ R=-aS-b\id_{E},&\label{relation1.4comp}\\
&g(PX,Y)=-g(X,PY)-ag(X,Y),&\label{relation1.5comp}\\
&g(QX,\xi)=-g(X,R\xi),&\label{relation1.6comp}\\
&g(S\xi,\nu)=-g(\xi,S\nu)-ag(X,Y),&\label{relation1.7comp}\\
&\nabla_X(PY)-P(\nabla_XY)=A_{QY}X+R(B(X,Y)),& \label{relation2.1comp}\\
&\nabla^{\perp}_X(QY)-Q(\nabla_XY)=S(B(X,Y))-B(X,PY),& \label{relation2.2comp}\\
&\nabla_X^{\perp}(S\nu)-S(\nabla^{\perp}_X\nu)=-B(R\nu,X)-Q(A_{\nu}X),& \label{relation2.3comp}\\
&\nabla_X(R\nu)-R(\nabla^{\perp}_X\nu)=-P(A_{\nu}X)+A_{S\nu}X.&\label{relation2.4comp}
\end{align}
{\it Proof:} 
Writing $J$ as a matrix by blocks with respect to the decomposition $T\widetilde{M}=TM\oplus E$, we have
$$J=\left( 
\begin{array}{cc}
P&R\\
Q&S
\end{array}\right)$$
and so 
$$J^2=\left( 
\begin{array}{cc}
P^2+R\circ S&P\circ R+ R\circ S\\ \\
Q\circ P+S\circ Q& Q\circ R+S^2
\end{array}
\right).$$
\end{prop}
The identities \eqref{relation1.1comp}-\eqref{relation1.4comp} are immediate from this and the relation $J^2+aJ+b\id=0$.\\
Moreover, the relations \eqref{relation1.5comp}-\eqref{relation1.7comp} come directly from the fact that for any $X,Y\in T\widetilde{M}$, $g(JX,Y)=g(X,JY)$.\\
Finally, \eqref{relation2.1comp}-\eqref{relation2.4comp} are consequences of the fact that $J$ is parallel and this is strictly the same as in Proposition \ref{allrelations}.\hfill$\square$\\ \\
From this proposition, we deduce immediately the following identities.
\begin{cor}\label{propgab}
For any $X,Y\in\Gamma(TM)$ and any $\nu,\xi\in\Gamma(E)$, we have
\begin{align}
&g(PX,PY)+g(QX,QY)=bg(X,Y),&\label{relationcor1}\\
&g(S\xi,S\nu)+g(R\xi,R\nu)=bg(\xi,\nu),&\label{relationcor2}\\
&g(PX,R\xi)+g(QX,S\xi)=0.&\label{relationcor3}
\end{align}
\end{cor}

\hfill$\square$\\ \\
We finish this section by considering two particular cases, namely the hypersurfaces and the invariant submanifolds. First, we consider invariant submanifolds.
\begin{definition}
A submanifold $M$ into the Riemannian  $(\widetilde{M},g)$ with a parallel Riemannian $(a,b)$-complex metallic structure is called invariant with respect to $J$ if $J(T_xM)\subset T_xM$ for all $x\in M$.
\end{definition}
First, we give the following proposition coming from Proposition \ref{allrelationscomp} for invariant submanifolds.
\begin{prop}
If $M$ is an invariant submanifold with respect to $J$, then the operators $Q$ and $R$ vanish and the operators $P$ and $S$ satisfy the following equations for all $X,Y\in TM$ and all $\xi,\nu\in E$.
\begin{align}
&P^2=-aP-b\id_{TM},&\label{relation1.1invcomp}\\
&S^2=-aS-bq\id_{E},&\label{relation1.4invcomp}\\
&g(PX,Y)=-g(X,PY)-ag(X,Y),&\label{relation1.5invcomp}\\
&g(S\xi,\nu)=-g(\xi,S\nu)-ag(\xi,\nu),&\label{relation1.7invcomp}\\
&\nabla_X(PY)-P(\nabla_XY)=0,& \label{relation2.1inv}\\
&S(B(X,Y))=B(X,PY),& \label{relation2.2inv}\\
&\nabla_X^{\perp}(S\nu)-S(\nabla^{\perp}_X\nu)=0,& \label{relation2.3inv}\\
&P(A_{\nu}X)=A_{S\nu}X.&\label{relation2.4inv}
\end{align}
In particular, $P$ and $S$ are parallel Riemannian $(a,b)$-complex metallic structure respectively on $TM$ and $E$.
\end{prop}
{\it Proof:} The proof is immediate from Proposition \ref{allrelationscomp} with the fact that $Q=0$ and $R=0$.\hfill$\square$\\ \\ 
Now, we can proof the following propositions. First, we have a relation between invariant submanifolds and metallic structures.
\begin{prop}
Let $(M,g)$ be a Riemannian manifold isometrically immersed into a Riemannian manifold $(N,\widetilde{g})$ carrying a Riemannian $(a,b)$-complex metallic structure $\widetilde{P}$. Let $P,Q,R,S$ the four operators induced on $M$ by $\widetilde{P}$. Then, $M$ is invariant if and only if $P$ is a non-trivial Riemannian $(a,b)$-complex metallic structure on $M$.\\
In particular, invariant submanifolds of Riemannian manifolds with Riemannian $(a,b)$-complex metallic structure are even-dimensional.
\end{prop}
{\it Proof:} From \eqref{relation1.1comp}, we have $P^2+R\circ Q=-aP-b\id_{TM}$. Hence, if $M$ is invariant, $Q=0$ and so $P^2=-aP-b\id_{TM}$, that is, $P$ is a $(a,b)$-complex metallic structure. Conversely, if $P$ is a $(a,b)$-complex metallic structure, then, we have $R\circ Q=0$. But since $-R$ is the dual of $Q$, this implies immediately that $Q=0$ and so $M$ is invariant.\\
The fact that an invariant submanifold is even dimensional is a direct consequence of the first part together with Proposition \ref{equivstructurecomp}. \hfill$\square$.

\begin{prop}
Any invariant submanifold of a Riemannian manifold with parallel Riemannian $(a,b)$-complex metallic structure is minimal.
\end{prop}
{\it Proof:} We consider the map $j:TM\lgra TM$ defined by $jX=\frac{2}{\delta}PX+\frac{a}{\delta}X$. From Corollary \ref{propgab}, we deduce the following elementary lemma:
\begin{lemma}\label{lembaseon}
For any $X,Y\in TM$, we have $g(jX,X)=0$ and $g(jX,jY)=g(X,Y)$.
\end{lemma}
{\it Proof:}
First, we have 
\beQ
g(jX,X)&=&g\left( \frac{2}{\delta}PX+\frac{a}{\delta}X, X\right)\\
&=&\frac{2}{\delta}g(PX,X)+\frac{a}{\delta}g(X,X)\\
&=&-\frac{2}{\delta}g(X,PX)-\frac{2a}{\delta}g(X,X)+\frac{a}{\delta}g(X,X)\\
&=&-g(X,jX)=-g(jX,X).
\eeQ
Hence, we get $g(jX,X)=0$. Moreover, we also have
\beQ
g(jX,jY)&=&g\left( \frac{2}{\delta}PX+\frac{a}{\delta}X, \frac{2}{\delta}PY+\frac{a}{\delta}Y\right)\\
&=&\frac{4}{\delta^2}g(PX,PY)+\frac{a^2}{\delta^2}g(X,Y)+\frac{2a}{\delta^2}(g(PX,Y)+g(X,PY))\\
&=&\frac{4b}{\delta^2}g(X,Y)+\frac{a^2}{\delta^2}g(X,Y)-\frac{2a}{\delta^2}g(X,Y)\\
&=&\frac{4b-a^2}{\delta^2}g(X,Y)\\
&=&g(X,Y),
\eeQ
where we have used \eqref{relation1.5comp}, \eqref{relationcor1} and the fact that $\delta^2=b-4a^2$. This concludes the proof of the lemma.\hfill$\square$\\ \\
Moreover, we have the following fact coming from \eqref{relation2.2comp}. Since $M$ is invariant, we have for any $X,Y\in TM$,
\beqt\label{relationBPS}
S(B(X,Y))=B(X,PY).
\eeqt

From Lemma \ref{lembaseon}, we deduce the existence of an adapted local orthonormal frame $\{e_1,\cdots,e_{2k}\}$ of $TM$ such that for any $i\in\{1,\cdots,k\}$, $e_{2k}=je_{2k-1}$. Hence, using \eqref{relationBPS}, we have
\beQ
2kH&=&\sum_{j=1}^{2k}B(e_i,e_i)\\
&=&\sum_{i=1}^k\left(  B(e_{2i-1},e_{2i-1})+B(je_{2i-1},je_{2i-1})\right)\\
&=&\sum_{i=1}^k\left(  B(e_{2i-1},e_{2i-1})+B\left(\frac{2}{\delta}Pe_{2i-1}+\frac{a}{\delta}e_{2i-1},\frac{2}{\delta}Pe_{2i-1}+\frac{a}{\delta}e_{2i-1}\right)\right)\\
&=&\sum_{i=1}^k\left(  B(e_{2i-1},e_{2i-1})+\frac{4}{\delta^2}B(Pe_{2i-1},Pe_{2i-1})+\frac{4a}{\delta^2}B(Pe_{2i-1},e_{2i-1})+\frac{a^2}{\delta^2}B(e_{2i-1},e_{2i-1})\right)\\
&=&\sum_{i=1}^k\left(  B(e_{2i-1},e_{2i-1})+\frac{4}{\delta^2}S^2(B(e_{2i-1},e_{2i-1}))+\frac{4a}{\delta^2}S(B(e_{2i-1},e_{2i-1}))+\frac{a^2}{\delta^2}B(e_{2i-1},e_{2i-1})\right)\\
&=&\sum_{i=1}^k\left(  \frac{4}{\delta^2} (S^2+aS+b\id)(B(e_{2i-1},e_{2i-1}))+\left(1-\frac{4b-a^2}{\delta^2}\right)B(e_{2i-1},e_{2i-1})\right)\\
&=&0,
\eeQ
since $S^2+aS+b\id=0$ and $1-\frac{4b-a^2}{\delta^2}=0$. Hence, $M$ is minimal.\hfill$\square$

Now, we consider hypersurfaces. In this case, it is more convenient to consider the real-valued second fundamental form by taking the scalar product with the unit normal $\nu$. Therefore, the metallic structure $J$ on $\widetilde{M}$ implies the existence of a field of symmetric operators $P:TM\longrightarrow TM$, an vector field $V\in\Gamma(TM)$ and a smooth function $f$ on $M$. Note that $V$ and $f$ correspond to the tensors $R$ and $S$, respectively, in this case. The tensor $Q$ is just the dual $1$-from associated to $V$. These three objects satisfy the following relations.
\begin{prop}\label{allrelationshypcomp}
If $M$ is a hypersurface of a Riemannian manifold with parallel Riemannian $(a,b)$-complex metallic structure, then $P$, $V$ and $f$ satisfy for all $X,Y\in TM$,
\begin{align}
&P^2-\langle V,\cdot\rangle V=-P-b\id_{TM},&\label{relation1.1hypcomp}\\
&PV=-\frac{a}{2}V,&\label{relation1.2hypcomp}\\
&\|V\|^2=\frac{\delta}{2},&\label{relation1.3hypcomp}\\
&f=-\frac{a}{2},&\label{relation1.4hypcomp}\\
&g(PX,Y)=-g(X,PY)-ag(X,Y),&\label{relation1.5hypcomp}\\
&\nabla_X(PY)-P(\nabla_XY)=-\langle V,Y\rangle AX+\langle AX,Y\rangle V,& \label{relation2.1hypcomp}\\
&\nabla_XV=-P(AX)-\frac{a}{2}AX.& \label{relation2.3hypcomp}
\end{align}
\end{prop}
{\it Proof:} The proof is a direct consequence of Proposition \ref{allrelationscomp} with $R=V$ and $Q=-\langle V,\cdot\rangle$ with the fact that \eqref{relation1.7comp} gives directly $f=-\frac{a}{2}$.
\hfill$\square$\\ \\
From this proposition, we deduce the following result.
\begin{prop}\label{condmingeod}
If $(M,g)$ is a hypersurface of a Riemannian manifold with parallel Riemannian $(a,b)$-complex metallic structure, then
\begin{enumerate}
\item 
If $M$ is totally geodesic then $\nabla P=0$ and $\nabla V=0$.
\item If $\nabla P=0$ and $AV=0$, then $M$ is totally umbilical.
\end{enumerate}
\end{prop}
{\it Proof:} It is clear form \eqref{relation2.1hypcomp} and \eqref{relation2.3hypcomp} that if $M$ is totally geodesic, then $P$ and $V$ are parallel. Conversely, assume that $P$ and $AV=0$. Then we deduce from \eqref{relation2.1hypcomp} that for any $X\in\Gamma(TM)$,
$$AX=\frac{\langle AX,V\rangle V}{\|V\|^2}=0.$$
This concludes the proof.
\hfill$\square$

\subsection{Fundamental theorem of submanifolds in complex space forms}
Now, we consider submanifolds of complex space forms $\mathbb{M}_{\CC}(4c)$ of constant sectional curvature $4c$. First, we recall that the curvature of the complex space form $\mathbb{M}_{\CC}(4c)$ is given by

\begin{eqnarray}\label{curvcomp1}
\widetilde{R}(X,Y)Z&=&c\Bigg[\left\langle Y,Z\right\rangle X-\left\langle X,Z\right\rangle  Y+\left\langle JY,Z \right\rangle JX -\left\langle JX,Z\right\rangle  JY,+2\left\langle X,JY\right\rangle JZ,\Bigg].
\end{eqnarray}
where $J$ is the complex structure of $\mathbb{M}_{\CC}(4c)$. Equivalentely, the curvature can be expressed in terms of $(a,b)$-complex metallic structure. Namley, we have
\begin{eqnarray}\label{curvcomp2}
\widetilde{R}(X,Y)Z&=&c\left(1+\frac{a^2}{\delta^2}\right)(\langle Y,Z\rangle X-\langle X,Z\rangle Y) +\frac{2ac}{\delta^2}(\langle Y,Z\rangle \widetilde{P}X-\langle X,Z\rangle \widetilde{P}Y)\nonumber\\
&&+\frac{2ac}{\delta^2}(\langle  \widetilde{P}Y,Z\rangle X-\langle  \widetilde{P}X,Z\rangle Y)+\frac{4c}{\delta^2}(\langle \widetilde{P}Y,Z\rangle \widetilde{P}X-\langle \widetilde{P}X,Z\rangle \widetilde{P}Y)\nonumber\\
&&+\frac{2ca^2}{\delta^2}\langle X,Y\rangle Z +\frac{4ac}{\delta^2}(\langle X,Y\rangle \widetilde{P}Z+\langle X,\widetilde{P}Y\rangle Z)+\frac{8c}{\delta^2}\langle X,\widetilde{P}Y\rangle \widetilde{P}Z
\end{eqnarray}
where $\widetilde{P}$ is the canonical $(a,b)$-complex metallic structures associated with $J$ and by $\widetilde{P}=\frac{\delta}{2}J-\frac{a}{2}\id$. Hence, the Gauss, Ricci and Codazzi equations are given by
\beq\label{Gausscomp}
R(X,Y)Z&=&c\left(1+\frac{a^2}{\delta^2}\right)(\langle Y,Z\rangle X-\langle X,Z\rangle Y) +\frac{2ac}{\delta^2}(\langle Y,Z\rangle PX-\langle X,Z\rangle PY)\nonumber\\
&&+\frac{2ac}{\delta^2}(\langle  PY,Z\rangle X-\langle  PX,Z\rangle Y)+\frac{4c}{\delta^2}(\langle PY,Z\rangle PX-\langle PX,Z\rangle PY)\nonumber\\
&&+\frac{2ca^2}{\delta^2}\langle X,Y\rangle Z +\frac{4ac}{\delta^2}(\langle X,Y\rangle PZ+\langle X,PY\rangle Z)+\frac{8c}{\delta^2}\langle X,PY\rangle PZ\nonumber\\
&&+A_{B(Y,Z)}X-A_{B(X,Z)}Y,
\eeq

\beq\label{Codazzicomp}
(\nabla_XB)(Y,Z)-(\nabla_YB)(X,Z)&=&\frac{2ac}{\delta^2}(\langle Y,Z\rangle QX-\langle X,Z\rangle QY)+\frac{4c}{\delta^2}(\langle PY,Z\rangle QX-\langle PX,Z\rangle QY)\nonumber\\
&&+\frac{4ac}{\delta^2}\langle X,Y\rangle QZ+\frac{8c}{\delta^2}\langle X,PY\rangle QZ,
\eeq

\beq\label{Riccicomp}
R^{\perp}(X,Y)\nu&=&\frac{4c}{\delta^2}(\langle QY,\nu\rangle QX-\langle QX,\nu\rangle QY)+\frac{2ca^2}{\delta^2}\langle X,Y\rangle \nu \nonumber\\
&&+\frac{4ac}{\delta^2}(\langle X,Y\rangle S\nu+\langle X,PY\rangle \nu)+\frac{8c}{\delta^2}\langle X,PY\rangle S\nu\nonumber\\
&&+B(A_{\nu}Y,Z)-B(A_{\nu}X,Z).
\eeq

Now, we can state the compatibility equations for isometric immersion into complex space forms. Namely, we have:
\begin{definition}\label{defcompCPn}
We say that $(M,g,E,g^E,\nabla^E,B,P,Q,R,S)$ satisfies the compatiblity equations associated with $\mathbb{M}^N_{\CC}(4c)$ if 
\begin{enumerate}

\item equations \eqref{relation1.1comp}-\eqref{relation2.4comp} are satisfied,
\item the Gauss, Codazzi and Ricci equations \eqref{Gausscomp}, \eqref{Codazzicomp}, \eqref{Riccicomp} are satisfied.
\end{enumerate}
\end{definition}
Now, we have the following:
\begin{thm}\label{thm1comp}
Let $(M^n,g)$ be a simply connected Riemannian manifold and $E$ a $m$-dimensional vector bundle over $M$ endowed with a metric $g^E$ and a compatible conection $\nabla^E$ so that $M=N$ is even. Moreover, let $B:TM\times TM\longrightarrow E$ be a $(2,1)$-symmetric tensor and $P: TM \longrightarrow TM$, $Q: TM \longrightarrow E$, $R:E \longrightarrow TM$ and $S:E \longrightarrow E$ are four operators. If $(M,g,E,g^E,\nabla^E,B,P,Q,R,S)$ satisfies the compatiblity equations associated with $\mathbb{M}^{m+n}_{\CC}(4c)$ then, there exists an isometric immersion $\varphi:M\longrightarrow \mathbb{M}^{m+n}_{\CC}(4c)$ such that the normal bundle of $M$ for this immersion is isometric to $E$ and so that the second fundamental form $II$ and the normal connexion $\nabla^{\perp}$ are given by $B$ and $\nabla^E$. Precisely, there exists a vector bundle isometry $\widetilde{\varphi}: E\longrightarrow T^{\perp}\varphi(M)$ so that
$$II=\widetilde{\varphi}\circ B,$$
$$\nabla^{\perp}\widetilde{\varphi}=\widetilde{\varphi}\nabla^E.$$
Moreover, we have
$$\widetilde{P}(\varphi_*X)=\varphi_*(PX)+\widetilde{\varphi}(QX),$$
$$\widetilde{P}(\widetilde{\varphi}\nu)=\varphi_*(R\nu)+\widetilde{\varphi}(S\nu),$$
where $\widetilde{P}$ is the canonical $(a,b)$-complex metallic structure of $\mathbb{M}^{m+n}_{\CC}(c)$. Moreover, this immersion is unique up to an isometry of $\mathbb{M}^{m+n}_{\CC}(c)$.
\end{thm}
{\it Proof:} The proof is similar to the proof of Theorem \ref{thm1}. From $P,Q,R$ and $S$, we define the following four operators:
$$\left\{
\begin{array}{l}
j=\dfrac{2}{\delta}P+\dfrac{a}{2}\id_{TM},\\
h=\dfrac{2}{\delta}Q,\\
s=\dfrac{2}{\delta}R,\\
t=\dfrac{2}{\delta}S+\dfrac{a}{2}\id_{E}.
\end{array}
\right.$$
By straightforward computations, we show that $(M,g,E,g^E,\nabla^E,j,h,s,t)$ satisfies the compatibility equations for an isometric immersion into $\mathbb{M}_{\CC}(c)$ such that the complex structure is given by $j$, $h$, $s$ and $t$ (see \cite {PT,NR}). As in the proof of Theorem \ref{thm1}, we see easily that $P$, $Q$, $R$ and $S$ are then the restriction of the canonical $(a,b)$-complex metallic structure $\widetilde{P}$.\hfill$\square$.\\ \\
As for product spaces with Theorem \ref{corM2R2}, we are able to prove a spinoral version in low dimension, but with some differences. First, we obtain results for $\CC P^2$ which is not spin but only ${\rm Spin}^c$. Moreover, as precised in \cite{NR}, spinorial results are given only in complex and Lagrangian surfaces. Since, there is equivalence between complex immersion in the complex structure of $\CC P^2$ and invariant immersions for the canonical $(a,b)$-complex metallic structure, we can prove the following.
\begin{thm}\label{corCP2}
Let $(M^2,g)$ be an oriented Riemannian surface and $E$ an oriented vector bundle of rank $2$ over $M$ with scalar product $\langle\cdot,\cdot\rangle_E$ and compatible connection $\nabla^E$. We denote by $\Sigma=\Sigma M\otimes\Sigma E$ the twisted spinor bundle. Let $B:TM\times TM\lgra E$ be a bilinear symmetric map, $P:TM\lgra TM$ a $(a,b)$-complex metallic structure on $M$ and $S:E\lgra E$ a $(a,b)$-complex metallic structure on $E$. Assume that $S(B(X,Y))=B(X,P(Y))$ for all  $X\in \Gamma(TM)$ and consider $\{e_1, e_2\}$ an orthonormal frame of $TM$. Then, the  following two statements are equivalent.
\begin{enumerate}

\item There exists a ${\rm Spin}^c$  structure on $\Sigma M\otimes\Sigma E$ which auxiliary line bundle's curvature is given by $F^{M +E} (e_1, e_2):=F^M(e_1, e_2) + F^E(e_1, e_2) = 0$ and a spinor field $\varphi\in\Gamma(\Sigma M\otimes\Sigma E)$ satisfying for all $X \in \Gamma(TM)$,
\begin{eqnarray}\label{partspinorcomp}
\nabla_X\varphi&=&-\frac 12\zeta(X)\cdot\varphi-\frac{1}{2}X\cdot\varphi+\frac{i}{\delta}\left( \frac{a}{2}X-P(X)\right)\cdot\overline{\varphi},
\end{eqnarray}
such that $\varphi^+$ and $\varphi^-$ never vanish and where $\eta$ is given by
$$
\zeta(X)=\sum_{j=1}^2e_j\cdot B(e_j,X).
$$
\item There exists a local invariant isometric immersion of $(M^2,g)$ into $\CC P^2$ with $E$ as normal bundle and second fundamental form $B$ such that the canonical $(a,b)$-complex metallic structure of $\CC P^2$ over $M$ is given by $P$ and $S$ in the sense of Theorem \ref{thm1comp}.
\end{enumerate}
\end{thm}
{\it Proof:}
We set $j=\frac{2}{\delta}P+\frac{a}{\delta}\id_{TM}$ and $t=\frac{2}{\delta}S+\frac{a}{\delta}\id_{E}$. From Proposition \ref{equivstructurecomp} we have that $j$ and $t$ are complex structures if and only if $P$ and $S$ are parallel $(a,b)$-complex metallic structures (on $M$ and $E$ respectively). Moreover, a straightforward computation shows that $S(B(X,Y))=B(X,P(Y))$ if and only if $t(B(X,Y))=B(X,j(Y))$ and 
\eqref{partspinorcomp} is equivalent to 
$$\nabla_X\varphi=-\frac 12\zeta(X)\cdot\varphi-\frac{1}{2}X\cdot\varphi+\frac{i}{2}j(X)\cdot\overline{\varphi}.$$
Hence, applying Theorem 1.1 of \cite{NR} allows to get this theorem since being a complex immersion for the complex structure of $\CC P^2$ is equivalent to be invariant for the canonical $(a,b)$-complex metallic structure.
\begin{remark}
As for the product case, we want to point out that it is also possible to have spinorial charcaterizations for hypersurfaces into $\CC P^2$ but with the existence of two spinor fields (see \cite{NR2}). We do not write here the analogue for complex metallic structures for briefness.
\end{remark}
\subsection{Examples}
In this section, we give some examples of applications of Theorem \ref{thm1comp}. First, we recall briefly the description of the 3-dimensional homogeneous manifolds with 4-dimensional isometry group. Such a manifold is a Riemannian fibration over a simply connected 2-dimensional manifold with constant curvature $\kappa$ and such that the fibers are geodesic. We denote by $\tau$ the bundle curvature, which measures the default of the fibration to be a Riemannian product. In fact, $\tau$ can be identified to the O'Neill tensor which is a well-known skew-symmetric tensor defined on Riemannian submersion. When $\tau$ vanishes, we get a product manifold $\M^2(\kappa)\times\RR$. Here, we describe 3-homogeneous manifolds with 4-dimensional isometry group and  $\tau\neq0$. These manifolds are of three types: they have the isometry group of the Berger spheres if $\kappa>0$, of the Heisenberg group $Nil_3$ if $\kappa=0$ or of $\widetilde{PSL_2(\RR)}$ if $\kappa<0$. In the sequel, we denote these homoegenous manifolds by $\Ekt$. For further details, one can refer to \cite{Da2} for instance.
\\
\indent
Let $\Ekt$ be a 3-dimensional homogeneous manifold with 4-dimensional isometry group. \linebreak Assume that $\tau\neq0$, {\it i.e.}, $\Ekt$ is not a  product manifold $\mathbb{M}^2(\kappa)\times\RR$. As we said, $\Ekt$ is a Riemannian fibration over a simply connected 2-dimensional manifold with constant curvature $\kappa$ and such that the fibers are geodesic. Now, let $\xi$ be a unitary vector field tangent to the fibers. We call it the vertical vector field. This vector field is a Killing vector field (corresponding to the translations along the fibers).\\
\indent
We denote respectively by $\nabla$ and $R$ the Riemannian connection and the curvature tensor of $\Ekt$. The manifold $\Ekt$ admit a local direct orthonormal frame $\{e_1,e_2,e_3\}$ with $e_3=\xi$
and such that the Christoffel symbols $\overline{\Gamma}_{ij}^k=\left\langle \nabla_{e_i}e_j,e_k\right\rangle$ are
\beqt\label{christoffel}
\left\lbrace  
\begin{array}{l}
\overline{\Gamma}_{12}^3=\overline{\Gamma}_{23}^1=-\overline{\Gamma}_{21}^3=-\overline{\Gamma}_{13}^2=\tau,\\ \\
\overline{\Gamma}_{32}^1=-\overline{\Gamma}_{31}^2=\tau-\sigma, \\ \\
\overline{\Gamma}_{ii}^i=\overline{\Gamma}_{ij}^i=\overline{\Gamma}_{ji}^i=\overline{\Gamma}_{ii}^j=0,\quad\forall\,i,j\in\{1,2,3\},
\end{array}
\right. 
\eeqt
where $\sigma=\dfrac{\kappa}{2\tau}$. Then we have
$$[e_1,e_2]=2\tau e_3,\quad [e_2,e_3]=\sigma e_1,\quad [e_3,e_1]=\sigma e_2.$$
We will call $\{e_1,e_2,e_3\}$ the canonical frame of $\Ekt$. From (\ref{christoffel}), we see easily that for any vector field $X$,
\beqt
\nabla_Xe_3=\tau X\wedge e_3,
\eeqt
where $\wedge$ is the vector product in $\Ekt$, that is, for any $X,Y,Z\in\Gamma(T\,M)$,
$$\left\langle X\wedge Y,Z\right\rangle =\det_{\{e_1,e_2,e_3\}}(X,Y,Z).$$
Moreover, from \eqref{christoffel}, we deduce that the curvature tensor $R$ is given by
\beq\label{curvEkt}
 R(X,Y)Z&=&(\kappa-3\tau^2)\Big(\langle X,Z\rangle Y-\langle Y,Z\rangle X\Big)\nonumber\\
&&+(\kappa-4\tau^2)\Big(\langle Y,\xi\rangle\langle Z,\xi\rangle X+\langle Y,Z\rangle\langle X,\xi\rangle \xi -\langle X\xi\rangle\langle Z,\xi\rangle Y-\langle X,Z\rangle\langle Y,\xi\rangle\xi\Big).
\eeq
Moreover, $\Ekt$ is endowed with a Sasakian structure $(\phi, \xi,\eta)$, with $\phi=\nabla_{(\cdot)}\xi$ and $\eta=\xi^{\#}$. We define the operator $P:T\Ekt\lgra T\Ekt$ by $PX=\frac{\delta}{2}\phi X-\frac{a}{2}X$, the vector $V=\frac{\delta}{2}\xi$ and the function $f=-\frac{a}{2}$. Moreover, we set $A$ the $(1,1)$-tensor defined by $AX=\tau X+\frac{4\tau^2-\kappa}{\tau\delta^2}\langle X,V\rangle V.$ We will show that $(\Ekt,g,P,V,f)$ satisfies the compatibility equations for an isometric immersion into the complex space form $\mathbb{M}_{\CC}\left(\frac{\kappa-4\tau^2}{4}\right)$ of constant holomoprhic curvature $\kappa-4\tau^2$ with $(a,b)$-complex structure given by $P,V$, $f$ and $A$ as shape operator. First, we have the following trivial relations 
\begin{lemma}\label{lemEkt1}
$P$, $V$, $f$ and $A$ satisfy all the relations of Lemma \ref{allrelationshypcomp}.
\end{lemma}
{\it Proof:} The computations are straightforward from the definition of these four objects.\hfill$\square$\\ \\
We have this second lemma which  gives the Gauss equation.
\begin{lemma}\label{lemEkt2}
The curvature tensor $R$ of $\Ekt$ satisfies
\beQ
 R(X,Y)Z&=&c(1+\frac{a^2}{\delta^2})R_0(X,Y)Z+\frac{4c}{\delta^2}R_0(PX,PY)Z+\frac{2ac}{\delta^2}\Big(R_0(X,PY)Z+R_0(PX,Y)Z\Big)\\
 &&-2c\Big(\frac{4}{\delta^2}\langle X,PY\rangle  PZ+\frac{2a}{\delta^2}\langle X,PY\rangle Z+\frac{2a}{\delta^2}\langle X,Y\rangle PZ+\frac{a^2}{\delta^2}\langle X,Y\rangle Z\Big)\\
 &&+R_0(AX,AY)Z
\eeQ
with $c=\frac{\kappa}{4}-\tau^2$ and $R_0$ is the curvature tensor given by
 $$R_0(X,Y)Z=\langle X,Z\rangle Y-\langle Y,Z\rangle X.$$

 \end{lemma}
{\it Proof:} From \eqref{curvEkt}, we get easily that
$$\left\{\begin{array}{l}
R(e_1,e_2)e_1=(\kappa-3\tau^2)e_2\\
R(e_1,e_2)e_2=-(\kappa-3\tau^2)e_1\\
R(e_1,e_2)\xi=0\\
R(e_1,\xi)e_1=\tau\xi\\
R(e_1,\xi)e_2=0\\
R(e_1,\xi)\xi=-\tau e_1\\
R(e_2,\xi)e_1=0\\
R(e_2,\xi)e_2=\tau\xi\\
R(e_2.\xi)\xi=-\tau e_2.
\end{array}\right.$$
If we denote $R_1(X,Y)Z$ the right hand side term in the statement of the lemma, we have by a straightforward computation
\beQ
R_1(X,Y)Z&=&c\Big(R_0(X,Y)Z+R_0(\phi X,\phi Y)Z-2\langle X,\phi Y\rangle\phi Z\Big)+R_0(AX,AY)Z.
\eeQ
Now, using the fact that $\left\{\begin{array}{l}
\phi e_1=e_2\\
 \phi e_2=-e_1\\
 \phi\xi=0\end{array}\right.$ and $\left\{\begin{array}{l}
Ae_1=\tau e_1\\
 Ae_2=\tau e_2\\
 A\xi=\left(2\tau-\frac{\kappa}{4\tau}\right)\xi,\end{array}\right.$\\ we obtain easily that $R_1=R$. For instance, we have
 \beQ
 R_1(e_1,e_2)e_1&=&c\Big(R_0(e_1,e_2)e_1+R_0(\phi e_1,\phi e_2)e_1-2\langle e_1,\phi e_2\rangle\phi e_1\Big)+R_0(Ae_1,Ae_2)e_1\\
 &=&c\Big(R_0(e_1,e_2)e_1-R_0( e_2, e_1)e_1+2\langle e_1,e_1\rangle e_2\Big)+\tau^2R_0(e_1,e_2)e_1\\
 &=&4ce_2+\tau^2 e_2\\
 &=&(k-4\tau^2)e_2+\tau^2e_2\\
 &=&(k-3\tau^2)e_2\\
 &=&R(e_1,e_2)e_1.
 \eeQ
The other equalities are in the same spirit and straightforward.\hfill$\square$\\ \\
Finally, we have this third Lemma which gives the Codazzi equation.
\begin{lemma}\label{lemEkt3}
The tensor $A$ satisfies
\beQ
d^{\nabla}A(X,Y)&=&\frac{\kappa-4\tau^2}{4}\left(\frac{4}{\delta^2}(\langle V,X\rangle PY-\langle V,Y\rangle PX)+\frac{2a}{\delta^2}(\langle V,X\rangle Y-\langle V,Y\rangle X)\right)\\
&&-\frac{\kappa-4\tau^2}{\delta^2}\Big(2\langle PX,Y\rangle V+a\langle X,Y\rangle V\Big).
\eeQ
\end{lemma}

{\it Proof:} For any $X,Y\in\Gamma(T\Ekt)$, we have
\beQ
\nabla_X(AY)&=&\nabla_X\left(  \tau Y+\frac{4\tau^2-\kappa}{\tau\delta^2}\langle Y,V\rangle V\right)\\
&=&\tau\nabla_XY+\frac{4\tau^2-\kappa}{\tau\delta^2}\Big(\langle \nabla_XY,V\rangle +\langle Y,\nabla_XV\rangle V+\langle Y,V\rangle \nabla_XV\Big)\\
&=&\tau\nabla_XY+\frac{4\tau^2-\kappa}{\tau\delta^2}\Big(\langle \nabla_XY,V\rangle-\tau \langle Y,PX\rangle V-\frac{a\tau}{2}\langle Y,X\rangle V-\tau\langle Y,V\rangle PX-\frac{a\tau}{2}\langle Y,V\rangle X\Big)
\eeQ
where we have used that $\nabla_X V=\tau (-PX+\frac{a}{2}X)$. Hence, from the definition of $\nabla_X(AY)=\nabla_X(AY)-\nabla_Y(AX)-A[X,Y]$ and since $\nabla$ is torsion-free, we get

\beQ
d^{\nabla}A(X,Y)&=&\frac{\kappa-4\tau^2}{4}\left(\frac{4}{\delta^2}(\langle V,X\rangle PY-\langle V,Y\rangle PX)+\frac{2a}{\delta^2}(\langle V,X\rangle Y-\langle V,Y\rangle X)\right)\\
&&-\frac{\kappa-4\tau^2}{\delta^2}\Big(\langle PX,Y\rangle V-\langle PY,X\rangle V\Big).
\eeQ
Moreover, from\eqref{relation1.5hypcomp}, we have $\langle PX,Y\rangle V-\langle PY,X\rangle V=2\langle PX,Y\rangle V+a\langle X,Y\rangle V$, which concludes the proof.
\hfill$\square$\\ \\
Thus, Lemmas \ref{lemEkt1}, \ref{lemEkt2} and \ref{lemEkt3} ensure that $\Ekt$ with the given objects $P$, $V$, $f$ and $A$ satisy the compatibility equations of defintion \ref{defcompCPn} and by Theorem \ref{thm1comp}, $\Ekt$ is isometrically immersed into the complex space form $\mathbb{M}_{\CC}\left(\frac{\kappa-4\tau^2}{4}\right)$ of constant holomoprhic curvature $\kappa-4\tau^2$ such that the canonical $(a,b)$-complex structure on $\mathbb{M}_{\CC}\left(\frac{\kappa-4\tau^2}{4}\right)$ is given over $\Ekt$ by $P,V$ and $f$ and such that $A$ is shape operator.

\end{document}